\documentclass{amsart}
\usepackage{amsmath,amssymb,epsfig}
\usepackage{color}
\usepackage{mathrsfs}
\usepackage{hyperref}
\newcommand{\arxiv}[1]{\href{http://www.arXiv.org/abs/#1}{arXiv:#1}}
\DeclareMathAlphabet{\mathbbold}{U}{bbold}{m}{n}
\newcommand{\zero}{\mathbbold{0}}
\newcommand{\unit}{\mathbbold{1}}
\newtheorem{theorem}{Theorem}[section]
\newtheorem{proposition}[theorem]{Proposition}
\newtheorem{lemma}[theorem]{Lemma}
\newtheorem{corollary}[theorem]{Corollary}
\theoremstyle{definition}
\newtheorem{definition}[theorem]{Definition}
\theoremstyle{remark}
\newtheorem{remark}[theorem]{Remark}
\newtheorem{example}[theorem]{Example}
\newcommand{\figperso}[2]{\begin{figure}[htb]\begin{center} #1\end{center}
#2\end{figure}}

\newcommand{\Cellmuh}{C^\ell_{\nu(\mu)}}
\newcommand{\iY}{\mathsf{Y}}

\newcommand{\ptyP}{\mathscr{P}}
\def\d+{\dot{+}}
\def\card#1{\# #1}
\newcommand{\cart}{\!\times\!}
\newcommand{\restr}[2]{#1\,\cap\,#2\cart #2}
\def\Nil(#1:#2){N\big[\begin{smallmatrix}#1\cr #2\end{smallmatrix}\big]}
\def\Anil(#1:#2){A\big[\begin{smallmatrix}#1\cr #2\end{smallmatrix}\big]}
\def\Vnil(#1:#2){V\big[\begin{smallmatrix}#1\cr #2\end{smallmatrix}\big]}
\def\SMALLMATRIX#1{\left[\begin{smallmatrix} #1 \end{smallmatrix}\right]}
\newcommand{\epi}{\operatorname{epi}}
\newcommand{\eps}{\epsilon}
\newcommand{\sC}{\mathcal{C}}

\newcommand{\corn}[1]{\mathsf{R} (#1)}
\newcommand{\cont}{\sC}

\newcommand{\dex}[1]{\mathsf{e}(#1)}
\newcommand{\dexo}{\mathsf{e}}
\newcommand{\Sym}{\mathfrak{S}}
\newcommand{\R}{\mathbb{R}}
\newcommand{\Rbar}{\overline{\R}}
\newcommand{\N}{\mathbb{N}}
\newcommand{\Z}{\mathbb{Z}}
\newcommand{\C}{\mathbb{C}}
\newcommand{\sA}{\mathcal{A}}
\newcommand{\sF}{\mathcal{F}}
\newcommand{\sP}{\mathcal{P}}
\newcommand{\sL}{\mathcal{L}}
\newcommand{\sQ}{\mathcal{Q}}
\newcommand{\sM}{\mathcal{M}}
\newcommand{\sV}{\mathcal{V}}
\newcommand{\sW}{\mathcal{W}}
\newcommand{\sY}{\mathcal{Y}}
\newcommand{\sZ}{\mathcal{Z}}
\newcommand{\weakm}{\prec^{\rm w}}
\newcommand{\mrm}[1]{\text{\rm #1}}
\newcommand{\sat}{{\mrm{Sat}}}
\newcommand{\val}{\operatorname{val}}
\newcommand{\sgn}{\operatorname{sgn}}
\newcommand{\rmin}{\R_{\min}}
\newcommand{\rminb}{\overline{\R}_{\min}}
\newcommand{\rhomin}{\rho_{\min}}
\newcommand{\eval}[1]{\widehat{#1}}
\newcommand{\divex}[1]{\overline{#1}}
\newcommand{\polfun}{\rmin\{\iY\}}
\newcommand{\ainf}{\hat{A}}
\newcommand{\dinf}{D}
\newcommand{\sch}{\operatorname{Schur}}
\newcommand{\vex}{\operatorname{vex}}
\newcommand{\perm}{\operatorname{per}}
\newcommand{\tr}{\operatorname{tr}}
\newcommand{\diag}{\operatorname{diag}}
\newcommand{\GC}{G^c}
\renewcommand{\emptyset}{\varnothing}
\newcommand{\bydef}{\stackrel{\text{\rm def}}{=}}
\newcommand{\new}[1]{{\em #1}\index{#1}}
\newcommand{\set}[2]{\{#1\mid\,#2\}}
\newcommand{\tend}{\mathop{\longrightarrow}}
\title[Min-plus methods in eigenvalue perturbation theory]{Min-plus methods in eigenvalue perturbation theory
and generalised Lidski\u\i-Vi\v{s}ik-Ljusternik theorem}
\author{Marianne Akian}
\date{%4 February, 2004. Revised version,
February 14, 2006
}
\address{Marianne Akian, INRIA-Rocquencourt, Domaine de Voluceau,
78153 Le Chesnay C\'edex, France}
\email{Marianne.Akian@inria.fr}
\author{Ravindra Bapat}
\address{Ravindra Bapat, Indian Statistical Institute, New Delhi, 110016,
India}
\email{rbb@isid1.isid.ac.in}
\author{St\'ephane Gaubert}
\address{St\'ephane Gaubert, INRIA-Rocquencourt, Domaine de Voluceau,
78153 Le Chesnay C\'edex, France}
\email{Stephane.Gaubert@inria.fr}

\begin{document}
\maketitle
\begin{abstract}
We extend the perturbation theory of Vi\v{s}ik, Ljusternik and  Lidski\u\i\ 
for eigenvalues of matrices, using methods of min-plus algebra.
We show that the asymptotics of the eigenvalues of 
a perturbed matrix is governed by certain discrete optimisation problems,
from which we derive new perturbation formul\ae,
extending the classical ones and solving cases which were singular in previous
approaches. Our results include general weak majorisation inequalities,
relating leading exponents of eigenvalues of perturbed
matrices and min-plus analogues of eigenvalues. 
\end{abstract} 
\sloppy
%\tableofcontents 
\section{Introduction}
Let $\sA_{\eps}$ denote a $n\times n$ matrix whose
entries,
which are continuous functions of a parameter $\eps>0$,
satisfy
\begin{align}
(\sA_\epsilon)_{ij}=a_{ij}\epsilon^{A_{ij}}
+o(\epsilon^{A_{ij}})
\label{eq-0}
\end{align}
when $\eps$ goes to $0$,
where $a_{ij}\in \C$, and $A_{ij}\in \R\cup\{+\infty\}$.
(When $A_{ij}=+\infty$, this means by convention that 
$(\sA_\epsilon)_{ij}$ is identically zero.)
The goal of this paper is to give
first order asymptotics 
\[\sL^i_\eps\sim 
\lambda_i \eps^{\Lambda_i}\enspace,
\]
with $\lambda_i\in \C\setminus\{0\}$ and $\Lambda_i\in \R$, for each of the eigenvalues $\sL^1_\eps,
\ldots, \sL^n_{\eps}$ of $\sA_{\eps}$.
%in some generic cases. 

Computing the asymptotics of spectral elements
is a central problem of perturbation theory,
see~\cite{kato} and~\cite{baumgartel}.
For instance, when the entries of $\sA_{\eps}$ 
have Taylor (or, more generally, Puiseux) series
expansions in $\eps$, 
the eigenvalues $\sL^i_{\eps}$ have Puiseux
series expansions in $\eps$, which can be computed by
applying the Newton-Puiseux algorithm to the
characteristic polynomial of $\sA_{\eps}$.
The leading exponents $\Lambda_i$ of the eigenvalues
of $\sA_{\eps}$ are
the slopes of the associated Newton polygon:
the difficulty is to determine these slopes 
from $\sA_{\eps}$.

The case of a linear perturbation of degree one
\[
\sA_\epsilon = \sA_0 + \epsilon b \enspace,
\qquad\sA_0,b\in \C^{n\times n}\enspace,
\]
has been particularly studied.
It suffices to consider the case where $\sA_0$ is nilpotent,
which is the object of a theory initiated
by Vi\v sik and Ljusternik~\cite{vishik}
and completed by Lidski\u\i~\cite{lidskii}.
Their result shows that
for generic values of the entries of $b$, the
exponents $\Lambda_i$ are the inverses of the dimensions
of the Jordan blocks of $\sA_0$. Then,
the coefficients $\lambda_i$ can be obtained from
the eigenvalues of certain Schur complements built from
the matrices $\sA_0$ and $b$. 

%The reader may consult the
%survey of Moro, Burke, and Overton~\cite{mbo97}
%for a complete account of this theory. 

However, the construction of Vi\v sik, Ljusternik and Lidski\u\i\
has many singular cases, in which the Schur complements
do not exist, and so, their approach does not apply to non-generic
situations, such as the case when
the matrix $b$ has a sparse or structured pattern.

The problem of generalising the theorem of
 Vi\v sik, Ljusternik and Lidski\u\i, i.e., of
``categorising all possible behaviours
as a function of the perturbation $b$'', to quote the introduction
of the article of Ma and Edelman~\cite{maedelman}, has received much attention.
Their article solves
cases where $\sA_0$ and $b$ have certain Jordan and Hessenberg structures,
respectively. This problem is also considered
in the survey of Moro, Burke, and Overton~\cite{mbo97},
which includes a slight refinement of Lidski\u\i's result
together with an extension in special cases.
Similar problems have been raised
for matrix pencils, see in particular Najman~\cite{najman}.
See also Edelman, Elmroth and K\aa gstr\"om~\cite{eek1,eek2}
for a geometric point of view. Numerical motivation can be found
there, as well as in the theory of pseudospectra, see Trefethen
and Embree~\cite{tref} for an overview. 
In~\cite{murota}, Murota gave an algorithm to compute the Puiseux
series expansions of the eigenvalues of a matrix whose entries
are given by polynomials (or even formal series) in some indeterminate.
This algorithm sheds some light on the problem raised by Ma and Edelman
(although this problem is not considered there).

In this paper, we use min-plus algebra to 
give elements of answer to the problem raised
by Ma and Edelman.

To describe our results,
let us recall that the min-plus semiring, $\rmin$,
is the set $\R\cup\{+\infty\}$, equipped
with the addition $(a,b)\mapsto \min(a,b)$
and the multiplication $(a,b)\mapsto a+b$.
Many of the classical algebraic constructions
have interesting min-plus analogues. In particular,
the characteristic polynomial function 
of a matrix $B\in \rmin^{n\times n}$,
already introduced by Cuninghame-Green~\cite{cuning83},
can be defined as the function which associates to a scalar
$x$ the permanent, in the min-plus sense,
of the matrix $ x I\oplus B$,
where $I$ is the min-plus identity matrix,
``$\oplus$'' denotes
the min-plus addition, and the concatenation denotes
the min-plus multiplication. 
The permanent, in the min-plus sense, of a matrix $B$,
is the value of an optimal
assignment in the weighted
bipartite graph associated with $B$.
A result of Cuninghame-Green
and Meijer~\cite{cuning80} shows that a min-plus polynomial
function $p(x)$ can be factored uniquely as 
$p(x)=a(x\oplus x_1)\cdots (x\oplus x_n)$,
where $a,x_1,\ldots,x_n\in \rmin$. The numbers $x_1,\ldots,x_n$, 
which coincide with the points of non-differentiability of $p$
(counted with appropriate multiplicities), are
called the {\em roots} or {\em corners} of $p$. 
The sequence of roots
of the min-plus characteristic polynomial of a
matrix $B\in\rmin^{n\times n}$
can be computed in polynomial time, by solving $O(n)$
optimal assignment problems, as shown
by Burkard and Butkovi\v{c}~\cite{bb02}.
The reader
seeking information on the min-plus semiring
may consult~\cite{cuning,maslov92,bcoq,maxplus94,cuning95b,guna96,maslovkolokoltsov95,maxplus97,pin95,gondran02,litmas}.

We assume that $\sA_\epsilon$
 is given by~\eqref{eq-0}. This allows one to handle the case
of a perturbed matrix $\sA_\epsilon=\sA_0+\epsilon b$, where
the matrix  $b$ is non-generic.

The first main result of the present paper, Theorem~\ref{th-gen=}, shows that
the sequence of leading exponents of the eigenvalues
of the matrix $\sA_\epsilon$ is weakly
(super) majorised by the sequence of roots of the 
min-plus characteristic polynomial of the matrix of leading
exponents of $\sA_\epsilon$,
and that the equality holds for generic values of the coefficients
$a_{ij}$.

The proof of Theorem~\ref{th-gen=} relies
on a variant of the Newton-Puiseux
theorem in which the data are only assumed to have first order
asymptotics, that we state as Theorem~\ref{th3}
in a way which illuminates the role of min-plus algebra.
We consider the branches $\sY(\epsilon)$
solutions of the equation $\sP(\epsilon,\sY(\epsilon))=0$,
where $\sP(\epsilon,Y)=\sum_{j=0}^{n} \sP_j(\eps) Y^j$
and the $\sP_j(\epsilon)$ are continuous functions,
such that $\sP_j(\epsilon)= p_j\epsilon^{P_j}+o(\epsilon^{P_j})$,
with $p_j\in \C$ and $P_j\in \R\cup\{+\infty\}$.
We characterise the cases where this information is enough to determine the
first order asymptotics of the branches 
$\sY_1(\epsilon),\ldots,\sY_n(\epsilon)$. Then,
the leading exponents of the branches
are precisely the roots of the min-plus polynomial
$P(Y)=\bigoplus_{j=0}^n P_jY^j$: the leading exponents
of the classical roots are the min-plus roots.
Note that in this case, by Legendre-Fenchel duality,
the roots of the min-plus polynomial $P(Y)$
are precisely the slopes of the Newton-Polygon
classically associated to $\sP (\epsilon,Y)$.
Note also that the generic equality in Theorem~\ref{th-gen=}
could be derived from Murota's combinatorial
relaxation technique~\cite{murota}, which uses
a parametric assignment problem, whose value
is exactly the min-plus characteristic
polynomial.

Theorem~\ref{th-gen=} determines the generic leading
exponents of the eigenvalues of $\sA_\epsilon$, but it does not
determine the coefficients $\lambda_i$.
To compute these coefficients, we
define, in terms of eigenvalues of min-plus Schur complements,
a sequence of {\em critical values} of $A$,
that we characterise as generalised circuit means.
We show that the sequence
of roots of the min-plus characteristic polynomial
of $A$ is weakly majorised by the sequence of 
critical values of $A$ (Theorem~\ref{th-gam-bet}), and we characterise
the equality case in terms of the existence of disjoint
circuit covers, or perfect matchings,
in certain graphs.

Our second main result, Theorem~\ref{th-1}, shows that,
in the equality case of Theorem~\ref{th-gam-bet},
the coefficients $\lambda_i$ can be obtained in terms
of eigenvalues of certain Schur complements constructed
from the matrix $a$.
The theorem of Vi\v sik, Ljusternik and Lidski\u\i\,
is a special case of this result
(Corollary~\ref{cor-l}).
We give in Section~\ref{sec-singular} examples
of singular cases which can be solved by Theorem~\ref{th-1}.
In the remaining singular cases, different methods
should be used, along the lines of~\cite{murota,abg04}.

We also prove an asymptotic result for eigenvectors,
Theorem~\ref{th-2}, which is analogous to Theorem~\ref{th-1}. However, 
the combinatorial characterisation of the cases where Theorem~\ref{th-2}
determines the generic asymptotics of all the entries
of eigenvectors is lacking, see Section~\ref{rk-formal}.
Note that even when the first order asymptotics of an eigenvalue
is determined, a detailed asymptotic information on $\sA_\epsilon$
may be needed to determine the first order asymptotics of the corresponding
eigenvector, as shown in our earlier work~\cite{ABG96} which concerns
the special case of the Perron eigenvector.

\if{
We determine the asymptotics of the eigenvalues by
means of certain discrete optimisation problems,
which are variations on the shortest path and optimal assignment
problems. This allows
us to generalise the theorem of
Vi\v sik, Ljusternik and Lidski\u\i.
Min-plus algebra is the fundamental tool to establish these
results. 

Our method relies on a fundamental tool: min-plus algebra.

We shall see that some of Murota's results are intimatelty related
to min-plus results

We show that,
for generic values of the coefficients $a_{ij}$,
the Newton polygon of the characteristic polynomial
of $\sA_\epsilon$ can be computed (in polynomial time)
as a function of the matrix of leading exponents of
$\sA_\epsilon$, $A=(A_{ij})$, using a min-plus analogue of the characteristic
polynomial. Then, the leading exponents $\Lambda_i$ of the eigenvalues
of $\sA_\epsilon$ are readily obtained.
We also show that
the leading coefficients $\lambda_i$ of the eigenvalues
of $\sA_\epsilon$ can be computed by Schur complement
formul\ae\ extending the ones of Lidski\u\i, when
the matrix of leading exponents, $A$,
satisfies a structural condition expressed in terms
of the existence of perfect matchings in certain digraphs.
In a companion paper~\cite{abg04}, we show that,
when this structural condition is not satisfied, 
the leading coefficients $\lambda_i$ can still
be computed (for generic values of the $a_{ij}$),
but one has to renounce to Schur complement formul\ae.
In the remaining singular cases,
the knowledge of the first
order asymptotics $(\sA_\epsilon)_{ij}\sim a_{ij}\epsilon^{A_{ij}}$
is not enough to determine the first order asymptotics
of the eigenvalues of $\sA_\epsilon$, higher order
expansions of the entries $(\sA_\epsilon)_{ij}$ are needed.}\fi

The present results provide a new illustration
of the role of min-plus algebra in asymptotic analysis,
which was recognised by Maslov~\cite[Ch.~VIII]{maslov73}.
He observed that WKB-type or large deviation type asymptotics lead
to limiting equations, like Hamilton-Jacobi equations,
satisfying some idempotent superposition principle.
So, min-plus algebra arises as the limit of a deformation
of usual algebra.
This observation is at the origin
of idempotent analysis~\cite{maslov92,maslovkolokoltsov95,litvinov00}.
It has been used by Dobrokhotov, Kolokoltsov, and
Maslov~\cite{maslov92a,maslovkolokoltsov95} to obtain
precise large deviation asymptotics concerning the Green
kernel and the first eigenvalues of a class of linear
partial differential equations,
with application to the Schr\"odinger equation.

The same deformation has been identified by Viro~\cite{viro},
in relation with the patchworking method he developed
for real algebraic curves.
It appears in several recent works in ``tropical geometry'',
in particular, by Mikhalkin~\cite{mikhalkin,mikhalkincras},
Forsberg, Passare, and Tsikh~\cite{forsbergpassaretsikh},
Passare and Rullgard~\cite{passare},
and Speyer and Sturmfels ~\cite{speyer}, 
following the introduction of am\oe bas of algebraic
varieties by Gelfand, Kapranov, and Zelevinsky~\cite{gelfand}.
In these works, the relation between Newton polytopes
and min-plus or ``tropical'' polynomials is apparent.
We use the same relation in the version of the Newton-Puiseux theorem
concerning first order asymptotics that we
stated as Theorem~\ref{th3}. 

Relations between max-plus algebra and asymptotic
problems have also appeared in other contexts.
Puhalski~\cite{tolya01} applied idempotent techniques
to large deviations theory.  Friedland~\cite{friedland}
observed that the max-plus eigenvalue can be
obtained as a limit of the Perron root.
Olsder and Roos~\cite{olsder}
and De Schutter and De Moor~\cite{deschutter} used asymptotics
theorems to derive certain max-plus algebraic identities.

Finally, we note that Theorem~\ref{th-1} was announced in~\cite{gaubert01c},
and that the role of min-plus roots in the Newton-Puiseux theorem
was mentioned in~\cite{autrans}.
\section{Preliminaries} 
In this section, we recall some classical facts
of min-plus algebra and show preliminary results.
See for instance~\cite{bcoq} for more details.

The \new{min-plus semiring}, $\rmin$, is the set
$\R\cup\{+\infty\}$ equipped with the addition
$(a,b)\mapsto a\oplus b=\min(a,b)$ and the multiplication
$(a,b)\mapsto a\otimes b=a+b$. We shall denote by $\zero=+\infty$
and $\unit=0$ the zero and unit elements of $\rmin$, respectively.
We shall use the familiar algebraic conventions, in
the min-plus context. For instance, if $A,B$ are
matrices of compatible dimensions
with entries in $\rmin$, $(AB)_{ij}=(A\otimes B)_{ij}
=\bigoplus_{k} A_{ik}B_{kj}
=\min_{k} (A_{ik}+B_{kj})$, $A^2=A\otimes A$, etc.
Moreover, if $x\in \rmin\setminus\{\zero\}$, then
$x^{-1}$ is the inverse of $x$ for the $\otimes$ law, 
that is $-x$, with the conventional notation.
We shall also denote by $\rminb$ the complete 
min-plus semiring, which is the set $\R\cup\{\pm\infty\}$ equipped,
as $\rmin$, with the $\min$ and $+$ laws, with the
convention  $+\infty+(-\infty)=-\infty+(+\infty)=+\infty$.

\subsection{Min-plus spectral theorem}\label{spectral-sec}
To any $n\times n$ matrix $A$ with entries
in a semiring $S$, we associate the directed
graph $G(A)$, which has nodes $1,\ldots,n$ and an arc $(i, j)$
if $A_{ij}\neq \zero$, where $\zero$ denotes the zero
element of $S$. We say that $A$ is \new{irreducible} if $G(A)$
is strongly connected. 

We next recall some results of min-plus spectral
theory: the min-plus version of the Perron-Frobenius
theorem has been discovered, rediscovered, 
precised or extended, by many 
authors~\cite{cuning,vorobyev67,romanovski,gondran77,cohen83,maslov92}.
Recent presentations can be found in~\cite{bcoq,cuning95b,maxplus97,bapat98,agw04}.
\begin{theorem}[Min-plus eigenvalue, see e.g.~{\cite[Th.~3.23]{bcoq}}]\label{min-plus-spectral1}
An irreducible matrix $A\in (\rmin)^{n\times n}$
has a unique eigenvalue:
\begin{equation}
\rhomin(A) =\bigoplus_{k=1}^{n} \bigoplus_{i_1,\ldots,i_k}
(A_{i_1i_2}\cdots A_{i_ki_1})^{\frac{1}{k}}
\label{e-1} \enspace .
\end{equation}
\end{theorem}
With the usual notations,~\eqref{e-1} can be rewritten as
\begin{align*}
 \rhomin(A) =\min_{1\leq k\leq n} \min_{i_1,\ldots,i_k}
\frac{A_{i_1i_2}+\cdots+A_{i_ki_1}}{k}\enspace .
\end{align*}
If $p=(i_0,i_1,\ldots , i_k)$
is a path of $G(A)$, we denote by $|p|_A=
A_{i_0i_1}+ \cdots + A_{i_{k-1}i_k}$ the \new{weight}
of $p$, and by $|p|=k$ its \new{length}.
Since any circuit of $G(A)$ can be decomposed in elementary circuits,
which are of length at most $n$, $\rhomin(A)$ is
the \new{minimal circuit mean}:
\begin{equation}
 \rhomin(A) =\min_{c\mrm{ circuit in } G(A)} \frac{|c|_A}{|c|}
\enspace .
\label{e-2-2} 
\end{equation}
We say that a circuit $c=(i_1,i_2,\ldots, i_k,i_1)$ of
$G(A)$ is \new{critical} if $c$  attains
the minimum in~\eqref{e-2-2}, and we call critical
the nodes and arcs of this circuit.
The critical nodes and critical arcs form the \new{critical graph},
$\GC(A)$. We call \new{critical classes} the strongly connected 
components of $\GC(A)$. We will also use the name ``critical
class'' for the set of nodes of a critical class.

The \new{Kleene's star} of a matrix $A\in \rmin^{n\times n}$
is defined by
\[
A^*=I \oplus A \oplus A^2 \oplus \cdots \in
\rminb^{n\times n} \enspace,
\]
i.e.\ $(A^*)_{ij}=\inf_{k\geq 0} (A^k)_{ij}$, where $I=A^0$ is the identity 
matrix (we shall use the same notation $I$ for the identity matrix
of $\rmin^{n\times n}$, and for the identity matrix of $\C^{n\times n}$, 
for any $n$).
\begin{proposition}[See e.g.{~\cite[Th.~3.20]{bcoq}}]
\label{prop-finite}
All the entries of $A^*$ are $>-\infty$ if, and only if, 
$\rhomin(A)\geq 0$. Moreover, when $\rhomin(A)\geq 0$,
\[
A^*=I\oplus A\oplus \cdots \oplus A^{n-1}\enspace .
\]
\end{proposition}

\begin{theorem}[Min-plus eigenvectors, see e.g.{~\cite[Th.~3.100]{bcoq}}]\label{min-plus-spectral2}
Let $A\in \rmin^{n\times n}$ be an irreducible matrix, and
consider $\widetilde{A}=\rhomin(A)^{-1}A$.
Any eigenvector of $A$ is a linear
combination of the columns $\widetilde{A}^*_{\cdot,j}$
corresponding to critical nodes $j$. 
More precisely, if we select (arbitrarily) one node $j$ per
critical class and take the corresponding column $\widetilde{A}^*_{\cdot,j}$,
we obtain a minimal generating set of the eigenspace of $A$.
\end{theorem}
(In Theorem~\ref{min-plus-spectral2}, and in the sequel, we
write $\widetilde{A}^*_{\cdot,j}$
the $j$-th column of $(\widetilde{A})^*$.)

Given a matrix $A\in \rmin^{n\times n}$ and a vector $V\in \rmin^n$,
we define the \new{saturation graph}, $\sat(A,V)$, 
which has nodes $1,\ldots,n$, and an arc $(i, j)$ 
if $(AV)_i=A_{ij}V_j$ (that is  $(AV)_i=A_{ij}+V_j$ with the usual
notations).
The following simple result relates the
critical graph and the saturation graph.
\begin{proposition}[See e.g.~{\cite[Th.~3.98]{bcoq}}]
\label{survect}
Let $A\in \rmin^{n\times n}$ be an irreducible matrix 
with eigenvalue $\alpha$, and let $V\in \rmin^{n}\setminus\{\zero\}$.
If $AV=\alpha V$, then the strongly connected components of $\sat(A,V)$ are
exactly the strongly connected components of $G^c(A)$.
\end{proposition}

In fact, Theorem~3.98 of ~\cite{bcoq} only shows
that any circuit of the saturation graph belongs
to the critical graph, but the converse is 
straightforward. 

The following elementary result is a special
version of a maximum principle for ergodic control
problems, see~\cite[Lemma~3.3]{spectral2} for more
background, and~\cite[Lemma~1.4]{gg} for a proof
in the min-plus case.
\begin{proposition}
\label{survect2}
Let $A\in \rmin^{n\times n}$ be an irreducible matrix 
with eigenvalue $\alpha$, and let $V\in \rmin^{n}$.
If $AV\geq \alpha V$, then $(AV)_i=\alpha V_i$ for all critical nodes
$i$ of $A$.
\end{proposition}

The saturation graphs associated to the generators
of the eigenspace have a remarkable structure.
Say that a strongly connected component
$C$ of a graph is \new{final} if for each node $i$,
there is a path from $i$ to $C$, and if there is no arc leaving
$C$.
\begin{proposition} \label{form-sat}
Let $A\in \rmin^{n\times n}$ be an  irreducible matrix 
with eigenvalue $\alpha$, let $\widetilde{A}=\alpha^{-1} A$,
let $C$ be a critical class of $A$, and let $V$ be an eigenvector
of $A$. The following assertions are equivalent:
\begin{enumerate}
\item\label{proportional} $V$ is proportional to $\widetilde{A}^*_{\cdot j}$, for some
$j\in C$;
\item\label{final} $C$ is the unique final class of $\sat(A,V)$.
\end{enumerate}
\end{proposition}
\begin{proof}
We first prove \ref{proportional}$\implies$\ref{final}.
It is enough to consider the case when $V=\widetilde{A}^*_{\cdot j}$.
Since $A$ is irreducible, 
all the entries of $\widetilde{A}^*$ 
are $<+\infty$. 
Moreover, since $\rhomin(\widetilde{A})=0$, 
Proposition~\ref{prop-finite} yields
$\widetilde{A}^*=I\oplus \widetilde{A}\oplus \cdots \oplus 
\widetilde{A}^{n-1}$.
Hence, for all $i\neq j$, there exists a
path $p=(i_0=i,i_1,\ldots, i_k=j)$ from $i$ to $j$,
with length $1\leq k\leq n-1$,
and minimal weight, that is $\widetilde{A}^*_{ij}=\widetilde{A}_{i_0 i_1}
\cdots  \widetilde{A}_{i_{k-1} i_k}$.
By Bellman's optimality principle,
for all $0\leq l\leq m\leq k$,
the sub-path $(i_l,\ldots,i_m)$
has minimal weight: $\widetilde{A}^*_{i_l i_m}=\widetilde{A}_{i_l i_{l+1}}
\cdots \widetilde{A}_{i_{m-1} i_m}$.
Then, $\widetilde{A}^*_{i_l j}=\widetilde{A}_{i_l i_{l+1}}
\widetilde{A}^*_{i_{l+1} j}$, that is,
$\alpha V_{i_l}=A_{i_l i_{l+1}}V_{i_{l+1}}$, and 
$(i_l,i_{l+1})\in \sat(A, V)$ for all $l=0,\ldots,k-1$.
So for $i\neq j$, there is a path from $i$ to $j$ in $\sat(A, V)$.

Assume, by contradiction, that there exists $k\in C$ and $l\not\in C$
such that $(k,l)\in \sat(A, V)$.
Since $l\neq j$, there is a path from $l$ to $j$ in $\sat(A, V)$,
and since $C$ is a strongly connected component of $\sat(A, V)$
(by Proposition~\ref{survect}), there is a path from $j$ to $k$
in $\sat(A, V)$, which yields a circuit of $\sat(A, V)$
passing through $C$ and $k\not\in C$.
This contradicts the fact that $C$ is a
strongly connected component of $\sat(A, V)$.

We finally prove~\ref{final}$\implies$\ref{proportional}.
Assume that $C$ is the unique final class of $\sat(A,V)$,
and let us fix $j\in C$. 
Then, for each $i$, we can find a path $(i_0=i,\ldots,i_k=j)$
from $i$ to $j$ in $\sat(A,V)$, so that
$V_{i_0}=\widetilde{A}^*_{i_0i_1}V_{i_1}$,
\ldots,
$V_{i_{k-1}}=\widetilde{A}^*_{i_{k-1}i_k}V_{i_k}$.
Hence, $V_i= \widetilde{A}^*_{i_0i_1}
\cdots \widetilde{A}^*_{i_{k-1}i_k}V_{j}
\leq \widetilde{A}^*_{ij} V_j$. The other inequality
holds, since $V=\widetilde{A}V$ implies 
$V=\widetilde{A}^*V$. Thus, 
$V=  \widetilde{A}^*_{\cdot j} V_j$
is proportional to $\widetilde{A}^*_{\cdot j} $.
\end{proof}

\subsection{Min-plus polynomials}\label{sec-minpol}
We recall here some results about formal polynomials and 
polynomial functions over $\rmin$, and in particular 
a min-plus analogue of ``the fundamental theorem of algebra'',
which is due to Cuninghame-Green and Meijer~\cite{cuning80}.
The connection between the min-plus evaluation
morphism and the Fenchel transform, was already
observed in~\cite{cohen89b} and~\cite[Section~ 3.3.1]{bcoq}.

We denote by $\rmin[\iY]$ the semiring of formal polynomials with 
coefficients in $\rmin$ in the indeterminate $\iY$: a
\new{formal polynomial} $P\in \rmin[\iY]$ is nothing but
a sequence $(P_k)_{k\in\N}\in \rmin^\N$ such that $P_k=\zero$
for all but finitely many values of $k$. Formal
polynomials are equipped with the entry-wise
sum, $(P\oplus Q)_k=P_k\oplus Q_k$, 
and the Cauchy product, $(P Q)_k=\bigoplus_{0\leq i\leq k}P_i Q_{k-i}$.
As usual, we denote a formal polynomial
$P$ as a formal sum, $P=\bigoplus_{k=0}^{\infty} P_k \iY^k$.
We also define the \new{degree} and \new{valuation} of $P$:
$\deg P=\sup\set{k\in \N}{P_k\neq \zero}$,
$\val P=\inf\set{k\in \N}{P_k\neq \zero}$
($\deg P=-\infty$ and $\val P=+\infty$ if $P=\zero$).
To any $P\in \rmin[\iY]$, 
we associate the \new{polynomial
function} $\eval{P}:\rmin\to\rmin,\; y\mapsto \eval{P}(y)=
\bigoplus_{k=0}^{\infty} P_k y^k$, that is, with the
usual notation:
\begin{align}
\eval{P}(y)=\min_{k\in \N} (P_k + k y) \enspace .
\label{e-minpoly}
\end{align}
We denote by $\polfun$
the semiring of polynomial functions $\eval{P}$.
Contrary to the case of real or complex polynomials,
the evaluation morphism, $\rmin[\iY]\to\polfun,\; P\mapsto \eval{P}$ is
not injective. 
Indeed, the evaluation morphism is essentially a specialisation
of the Fenchel transform over 
$\R$:
\[ \sF : \Rbar^{\R} \to \Rbar^{\R},\; \sF (f)(y)=\sup_{x\in\R} (xy-f(x))
\enspace,
\]
since, for all $y\in\R$, 
$\eval{P}(y)= -\sF(P)(-y)$, where $P$ is extended to a function 
\begin{align}
\label{pasafunction}
P:\R\to \Rbar,\, x\mapsto P(x),
\mrm{ with } P(x) = \begin{cases}
P_k &\mrm{if }x=k\in \N\enspace ,\\
+\infty &\mrm{otherwise}
\end{cases}
\end{align}
It follows from~\eqref{e-minpoly}
that $\eval{P}$ is a concave nondecreasing function
with integer slopes.

In the sequel, we denote by $\vex f$ the convex
hull of a map $f:\R\to \Rbar$, and we
denote by $\divex{P}$ the formal polynomial
whose sequence of coefficients is obtained by
restricting to $\N$ the convex hull of the map $P:\R \to \Rbar$.
Thus, $\divex{P}_k=(\vex P)(k)$. The following
result is a special case of the Legendre-Fenchel
inversion theorem~\cite[Section~12]{rockafellar}.
\begin{proposition}\label{prop-mini}
If $P\in \rmin[\iY]$,
then $\divex{P}$ is the minimal formal polynomial $Q$ such that 
$\eval{Q}=\eval{P}$, we have
$\divex{\divex{P}}=\divex{P}$,
and $\divex{P}$ is given by
\[ \divex{P}_k= \sup_{y\in \R} (-k y +\eval{P}(y) ) \enspace .\]
\end{proposition}
\begin{theorem}[{\cite[Th.~3.43, 1 and 2]{bcoq}}]\label{th21}
A formal polynomial of degree $n$, $P \in \rmin[\iY]$, satisfies 
$P=\divex{P}$ if,
and only if, there exist $c_1\leq \cdots \leq c_n \in\rmin$ 
such that 
\[P = P_n (\iY\oplus c_1) \cdots (\iY\oplus c_n)\enspace
 .\]
The $c_i$ are unique and given,  by:
\begin{equation}\label{corners0}
c_i=\begin{cases}  P_{n-i}  (P_{n-i+1})^{-1}& \mrm{if } P_{n-i+1}\neq\zero\\
\zero &\mrm{otherwise,}
\end{cases}\qquad \mrm{for } i=1,\ldots, n\enspace .
\end{equation}
\end{theorem}
The min-plus analogue of the fundamental theorem of algebra due
to Cuninghame-Green and Meijer can be obtained
by applying Theorem~\ref{th21} to $\divex{P}$, since $\divex{\divex{P}}=\divex{P}$ and
$\eval{\divex{P}}=\eval{P}$.
\begin{theorem}[{\cite{cuning80}}]\label{th22}
Any polynomial function $\eval{P}\in \polfun$
can be factored in a unique way as
\begin{equation}
\label{th22.1} \eval{P}(y) = P_n (y\oplus c_1) \cdots 
(y\oplus c_n)\enspace ,
\end{equation}
with $c_1 \leq \cdots \leq c_n$.  
\end{theorem}

The $c_i$ are called the \new{roots} of $\eval{P}$. 
\typeout{CHECK cuning80}
(In~\cite{cuning80}, the term \new{corners} is used as a synonym of
root, we use the term of root which makes the analogy
with classical algebra clearer.)
The \new{multiplicity} of the root $c$ is the cardinality of the set
 $\set{j\in\{1,\ldots, n\}}{c_j = c}$. 
We shall denote by $\corn{\eval P}$ the sequence of roots: 
$\corn{\eval{P}}=(c_1,\ldots,c_n)$.
By extension, if $P\in\rmin[\iY]$ is a formal polynomial,
we will call \new{roots} of $P$ the roots of $\eval{P}$:
$\corn{P}=\corn{\eval{P}}$. By Proposition~\ref{prop-mini}, $\corn{P}=
\corn{\divex{P}}$.
Geometrically, the function $\divex{P}$ is the restriction to $\N$
of the convex function $\vex P$, which is
piecewise affine on its support, $[\val P,\deg P]$, and
$\eval{P}$ is concave, piecewise affine.
\begin{proposition}\label{prop-prenewton}
The roots $c\in\R$ of a formal polynomial $P\in \rmin[\iY]$
are exactly the points at which $\eval{P}$ is not differentiable.
They coincide with the opposites of the slopes of the affine parts
 of $\vex P:[\val P,\deg P]\to \R$.
The multiplicity of a root $c\in\R$
is equal to the variation of slope
of $\eval{P}$ at $c$, $\eval{P}'(c^-)-\eval{P}'(c^+)$,
and it coincides with the length of 
the interval where $\vex P$ has slope $-c$.
Moreover, $\zero$ is a root of $P$ if, and only if, $\eval{P}'(\zero^-):=
\lim_{c\to +\infty}\eval{P}'(c)\neq 0$.
In that case $\eval{P}'(\zero^-)$ is the multiplicity of $\zero$,
and it coincides with $\val P$.
\end{proposition}
\begin{proof}
The characterisation of the roots and of their multiplicities
in terms of $\eval{P}$ is due to Cuninghame-Green and Meijer~\cite{cuning80}.
It can be deduced from~\eqref{th22.1},
since when $c\in \R$, $\eval Q(y):=(y\oplus c)^k= k\min(y,c)$
has $c$ as unique point of non differentiability, with
$\eval{Q}'(c^-)=k$ and $\eval{Q}'(c^+)=0$. The case where $c=\zero$
is a straightforward consequence of~\eqref{th22.1}.
The characterisation of the roots and of their multiplicities
in terms of $\vex P$ follows from~\eqref{corners0},
since when $c_i\in \R$, $c_i= \divex{P}_{n-i}- \divex{P}_{n-i+1} 
=(\vex P)'(x)$  for all $x\in (n-i,n-i+1)$,
and $c_i=\zero\implies \divex{P}_{n-i}= \divex{P}_n c_1\cdots c_i=\zero$.
\end{proof}

The duality between roots and slopes in Proposition~\ref{prop-prenewton}
is a special case of the Legendre-Fenchel duality formula for
subdifferentials: $-c\in \partial (\vex P) (x)\Leftrightarrow
x\in \partial \sF (P) (-c) \Leftrightarrow
x\in \partial^+\eval{P} (c) $ where $\partial$ and $\partial^+$ denote
the subdifferential and superdifferential,
respectively~\cite[Th.~23.5]{rockafellar}.
\begin{lemma}\label{carac-ci}
Let $P=\bigoplus_{i=0}^n P_i \iY^i\in\rmin[\iY]$ be a formal polynomial
of degree $n$. Then, $\corn{P}=(c_1\leq \cdots \leq c_n)$ if, and only if,
$P\geq P_n (\iY\oplus c_1)\cdots (\iY\oplus c_n)$ and
\begin{equation}\label{carac-ci1}
 P_{n-i}=P_n c_1 \cdots c_i\quad \mrm{for all } i\in\{0,n\}\cup\set{i\in
\{1,\ldots , n-1\}}{c_i< c_{i+1}}\enspace .\end{equation}
In particular, $P_{n-i}=\overline{P}_{n-i}$ holds for
all $i$ as in~\eqref{carac-ci1}.
\end{lemma}
\begin{proof}
We first prove the ``only if'' part. 
If $\corn{P}=(c_1\leq \cdots \leq c_n)$,
then $\divex{P}=\divex{P}_n (\iY\oplus c_1)\cdots (\iY\oplus c_n)$ and
$\divex{P}_{n-i}= \divex{P}_n c_1\cdots c_i$ for all $i=1,\ldots n$.
Recall that $P$ defines a map $x\mapsto P(x)$ by~\eqref{pasafunction}.
By definition of $\vex P$,
the epigraph of $\vex P$, $\epi \vex P$,
is the convex hull of the epigraph of $P$,
$\epi P$. By a classical result~\cite[Cor~18.3.1]{rockafellar}, if $S$ is a 
set with convex hull $C$, any extreme point of $C$ belongs
to $S$. Let us apply this to $S=\epi P$ and
$C=\epi \vex P$.
Since $\divex{P}_{n-i}= \divex{P}_n c_1\cdots c_i$,
the piecewise affine map $\vex P$
changes its slope at any 
point $n-i$ such that 
$c_i<c_{i+1}$.
Thus, any point $(n-i,\vex{P}(n-i))$ with $c_i<c_{i+1}$
is an extreme point of $\epi \vex P$,
which implies that $(n-i,\vex{P}(n-i))\in \epi P$,
i.e., $P_{n-i}\leq \vex{P}(n-i)=\divex P_{n-i}$. Since
the other inequality is trivial by definition of the convex hull, we
have $P_{n-i}=\divex{P}_{n-i}$. 
Obviously, $P$ and $\divex P$ have the same
degree, which is equal to $n$, and they have the same
valuation, $k$. Then, $(n,\vex P(n))$
and $(k,\vex P(k))$ are extreme points
of $\epi \vex P$, and by the preceding argument,
$P_n=\divex P_n$, and $P_k=\divex P_k$.
Hence, $P_0=\divex P_0$, if $k=0$,
and $P_0=\divex P_0=+\infty$, if $k>0$.
We have shown ~\eqref{carac-ci1},
together with the last statement of the lemma.
Since $\divex{P}_n=P_n$ and $P\geq \divex{P}$, we also obtain 
$P\geq P_n (\iY\oplus c_1) \cdots (\iY\oplus c_n)$.

For the ``if'' part, assume that $P\geq P_n (\iY\oplus c_1) \cdots (\iY\oplus c_n)$
and that~\eqref{carac-ci1} holds. Since $Q=P_n (\iY\oplus c_1) \cdots
 (\iY\oplus c_n)$ is convex, and the convex hull map 
$P\mapsto \divex{P}$ is monotone, we must have $\divex{P}\geq \divex{Q}=Q$.
Hence, $ P \geq \divex{P}\geq Q$ and since $P_{n-i}=Q_{n-i}$ 
for all $i$ as in~\eqref{carac-ci1}, we must have $\divex{P}_{n-i}=
Q_{n-i}$, thus $\vex P(n-i)=\vex Q(n-i)$
at these $i$. Since $\vex P$ is convex, since $\vex Q$ is piecewise
affine and $\vex Q(j)=\vex{P}(j)$ for $j$ at the boundary of the domain
of $\vex Q$ and at all the $j$ where $\vex Q$ changes of slope, we must
have $\vex{P}=\vex Q$. Hence $\divex{P}=\divex{Q}=Q$
and $\corn{P}=\corn{\divex{P}}=\corn{Q}=(c_1,\ldots , c_n)$.
\end{proof}
The above notions are
illustrated in Figure~\ref{minpoly}, where
we consider the formal min-plus polynomial 
$P=\iY^3\oplus 5\iY^2 \oplus 6\iY \oplus 13$.
The map $j\mapsto P_j$, together with the 
map $\vex P$, are depicted at the left of the figure,
whereas the polynomial function $\eval{P}$ is depicted
at the right of the figure. We have $\divex{P}= \iY^3\oplus 3\iY^2 \oplus 6\iY\oplus 13
= (\iY\oplus 3)^2(\iY\oplus 7)$. Thus, the roots
of $P$ are $3$ and $7$, with respective multiplicities
$2$ and $1$. The roots are visualised at the right
of the figure, or alternatively, as the opposite
of the slopes of the two line segments at the
left of the figure. The multiplicities can be read
either on the map $\eval{P}$ at the right of the figure (the variation of slope
of $\eval{P}$ at points $3$ and $7$ is $2$ and $1$, respectively),
or on the map $\vex P$ at the left of the figure (as the respective horizontal
widths of the two segments). 
\figperso{\input{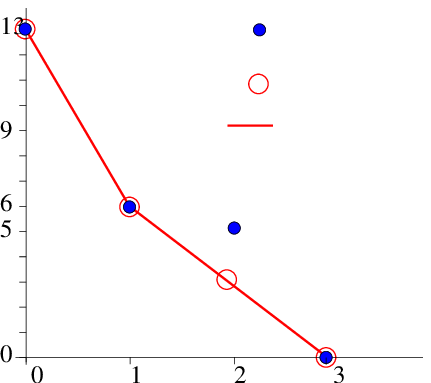}\hskip 2em \input{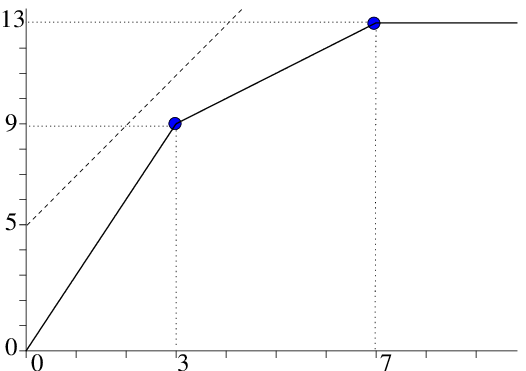}}{\caption{The formal min-plus polynomial $P=\iY^3\oplus 5\iY^2 \oplus 6\iY \oplus 13$ and its associated polynomial function $\eval P$.\label{minpoly}}}

\subsection{Schur complements}\label{sec-schur}
We recall here the definitions of conventional and min-plus
Schur complements.
We shall consider matrices indexed by ``abstract indices'':
if $L$ and $M$ are finite sets and  $S$ is a semiring, a
$L\times M$ matrix with values in $S$ is an element $A$ of
$S^{L\times M}$ and the entries of $A$ are denoted by $A_{ij}$ with $i\in L$
and $j\in M$. Moreover, for all $J\subset L$ and $K\subset M$, 
we denote by $A_{JK}$ the $J\times K$ submatrix of $A$:
$A_{JK}=(A_{jk})_{j\in J,\, k\in K}$. This definition applies to 
$n\times n$ matrices by taking $L=M=\{1,\ldots , n \}$.
Graphs of $L\times L$ matrices $A$ are defined as for $n\times n$ matrices
(see Section~\ref{spectral-sec}) with the only difference that the 
set of nodes is $L$.

\begin{definition}
Let $C\subset L$ be finite sets, and let $N=L\setminus C$.
If $a$ is a $L\times L$ matrix with entries in $\C$, 
and if $a_{CC}$ is invertible, the \new{Schur complement}
 of $C$ in $a$ is defined by
\[
\sch(C,a)= a_{NN}-a_{NC}(a_{CC})^{-1}a_{CN}
\enspace .
\]
\end{definition}

\begin{definition}
Let $C\subset L$ be finite sets, and let $N=L\setminus C$.
If $A$ is a $L\times L$ matrix  with entries in $\rmin$,
$\lambda\in \rmin\setminus\{\zero\}$,
and $\rhomin(\lambda^{-1}A_{CC})\geq 0$,
the min-plus $\lambda$-\new{Schur complement}
of $C$ in $A$ is defined by
\begin{align}
\sch(C,\lambda,A)=
 A_{NN} \oplus A_{NC}(\lambda^{-1}\!\!A_{CC})^*\lambda^{-1}
A_{CN}\enspace .
\label{e-def-mschur}
\end{align}
When $\lambda=\unit=0$,
we shall simply write $\sch(C,A)$ instead of
$\sch(C,\unit,A)$. 
\end{definition}
In fact, in the sequel, we shall mostly use min-plus
Schur complement corresponding to $\lambda=\rhomin(A)$.
The goal of the insertion of the normalising factors in~\eqref{e-def-mschur}
is to get the following homogeneity property:
\begin{align}\label{schur-homo}
\sch(C,\mu \lambda, \mu A) 
= \mu \sch(C,\lambda,A)\enspace, 
\end{align}
for all $\lambda,\mu\in \R$
such that $\lambda\leq \rhomin (A_{CC})$
and $\mu\lambda\leq \rhomin (A_{CC})$.

Using the same symbol, ``$\sch$'', both for
conventional and min-plus Schur complements
is not ambiguous: considering 
min-plus Schur complements of complex
matrices, or conventional Schur complements
of min-plus matrices, would be meaningless.

Both min-plus and conventional Schur
complements satisfy
\begin{equation}
\sch(C\cup C', a) = \sch(C,\sch(C',a))\label{e-ident}
\end{equation}
for all $L\times L$ matrices $a$, and
for all disjoint subsets of indices $C,C'\subset L$,
provided that the Schur complements are well defined
(if $\sch(C',a)$ is well defined, 
then the left hand side of \eqref{e-ident} exists if, and only if,
its right hand side exists). Of course,
~\eqref{e-ident} is a classical Gaussian elimination
identity, which is well known, both in conventional
algebra and in the min-plus algebra (the left hand side
and the right hand side of~\eqref{e-ident} are unambiguous
rational expressions, with elementary interpretations in terms
of paths, see for instance~\cite{lallement} for more background).

Finally, if $K\subset L$ and if $b$ is the $K\times K$ submatrix
of $a$, we shall sometimes write abusively $\sch(b,a)$, instead
of $\sch(K,a)$. 

We now give some graph interpretations of the weights and eigenvalues
of min-plus Schur complements.
Let $G$ be a graph with set of nodes $L$, let $C$ be a subset of $L$
and set $N=L\setminus C$. For all paths $p=(i_0,\ldots ,i_k)$ of $G$,
we denote by $|p|_C$ the number of
arcs of $p$ with initial node in $C$,
i.e., $|p|_C=\card{\set{0\leq m\leq k-1}{i_m\in C}}$,
where $\card{}$ denotes the cardinality of a set.
(All the path interpretations below have dual versions,
obtained by replacing ``initial'' by ``final''.)
We also denote by
$p\cap C$ the subsequence of $p$ obtained by deleting
the nodes not in $C$ ($p\cap C$ need not be a path of $G$). 
The following classical interpretation of Schur complements 
is an immediate consequence of the graph interpretation of the star.
\begin{lemma}\label{lem-schur}
Let $C\subset L$ be finite sets, and let $N=L\setminus C$.
Let $A$ be a $L\times L$ matrix  with entries in $\rmin$, and
$\lambda\in \rmin\setminus\{\zero\}$ be such that 
$\rhomin(A_{CC})\geq \lambda$.
Then, $p$ is a path in $G(\sch(C,\lambda,A))$ if, and only if,
there exists a path $p'$ in $G(A)$ with the same extremal nodes as $p$
and such that $ p'\cap N=p $.
Moreover, for all paths $p$ in $G(\sch(C,\lambda,A))$, we
have 
\[ |p|_{\sch(C,\lambda,A)}  = \min |p'|_A - \lambda |p'|_C  \enspace,
\]
where the minimum is taken over all the paths $p'$ of $G(A)$ that have 
the same extremal nodes as $p$ and satisfy  $ p'\cap N=p $.
In particular, $c$ is a circuit in $G(\sch(C,\lambda,A))$ if, and only if,
there exists a circuit $c'$ in $G(A)$ such that  $ c'\cap N=c $;
and for all circuits $c$ in $G(\sch(C,\lambda,A))$, we have 
\[ |c|_{\sch(C,\lambda,A)}  = \min |c'|_A - \lambda |c'|_C  \enspace,
\]
where the minimum is taken over all the circuits $c'$ of $G(A)$ such that
$c'\cap N=c$.
\end{lemma}
\begin{proposition}\label{prop1}
Let $C\subset L$ be finite sets, and let $N=L\setminus C$.
Let $A$ be a $L\times L$ matrix  with entries in $\rmin$, and
$\lambda\in \rmin\setminus\{\zero\}$ be such that 
$\rhomin(A_{CC})\geq \lambda$.
Then, 
\begin{align}
\rhomin(\sch(C,\lambda,A))=\min\frac{|c'|_A - \lambda |c'|_C}{
|c'| - |c'|_C}\label{eq-4}
\end{align}
where the minimum is taken over all the circuits $c'$ of $G(A)$ which are
not included in $C$.
Moreover, $c$ is a critical circuit of $\sch(C,\lambda,A)$ if, and only if,
there exists a circuit  $c'$ of  $G(A)$ such that $c'\cap N=c$
and $c'$ minimises~\eqref{eq-4}.
\end{proposition}
\begin{proof}
Using~\eqref{e-2-2} and Lemma~\ref{lem-schur}, we get
\begin{align*}
 \rhomin(\sch(C,\lambda,A))&=\min_{c\mrm{ circuit in }N}
\frac{|c|_{\sch(C,\lambda,A)}}{|c|}\\
&=\min_{c\mrm{ circuit in }N}\left(\min_{
c'\mrm{ circuit of } G(A),\,  c'\cap N=c} \frac{|c'|_A - \lambda |c'|_C}{|c|}
\right)\\
&= \min_{c'\mrm{ circuit of } G(A),\, c'\cap N\neq\emptyset} 
\frac{|c'|_A - \lambda |c'|_C}{|c'|-|c'|_C}\enspace ,
\end{align*}
since $|c'\cap N |=|c'|-|c'|_C $ for all circuits $c'$.
This yields~\eqref{eq-4}.
If $c$ is a critical circuit of $\sch(C,\lambda,A)$, then 
$\rhomin(\sch(C,\lambda,A))=(|c|_{\sch(C,\lambda,A)})/{|c|}$
and by Lemma~\ref{lem-schur},
there exists a circuit $c'$ of $G(A)$ such that $c'\cap N=c$ and
$|c|_{\sch(C,\lambda,A)}=|c'|_A-\lambda |c'|_C$. Since in that case,
$|c|=|c'\cap N |=|c'|-|c'|_C $, we deduce that $c'$ minimises~\eqref{eq-4}.
Conversely, if  $c'$ minimises~\eqref{eq-4}, then, $c=c'\cap N$ is
nonempty and by Lemma~\ref{lem-schur}, $c$ is a
circuit of $G(\sch(C,\lambda,A))$. Moreover, by Lemma~\ref{lem-schur}
again,
\[ \rhomin(\sch(C,\lambda,A))\leq 
\frac{|c|_{\sch(C,\lambda,A)}}{|c|}\leq 
\frac{|c'|_A - \lambda |c'|_C}{|c'|-|c'|_C}=\rhomin(\sch(C,\lambda,A))
\enspace , \]
thus $c$ is a critical circuit of $\sch(C,\lambda,A)$.
\end{proof}

Note that if $c'$ is a circuit in $C$,
that is if the denominator in~\eqref{eq-4} is zero, 
the numerator is necessarily nonnegative, since $\lambda\leq \rhomin(A_{CC})$.

\section{Min-plus polynomials, Newton-Puiseux theorem and 
generic exponents of eigenvalues}

\subsection{Preliminaries on exponents and general assumptions}
Let $\cont$ denote the set of continuous functions $f$ from some interval
$(0,\eps_0)$ to $\C$ with $\eps_0>0$,
such that $|f(\eps)|\leq \eps^{-k}$ on $(0,\eps_0)$,
for some positive constant $k$.
Since all the properties that we will prove in the sequel will hold
on some neighbourhoods of $0$, we shall rather use the ring of 
\new{germs} at $0$ of elements of $\cont$, which is obtained by quotienting
$\cont$ by the equivalence relation that identifies functions which coincide
on a neighbourhood of $0$.
This ring of germs will be also denoted by $\cont$.
For any germ $f\in\cont$, we shall abusively denote by $f(\eps)$ or $f_\eps$
the value at $\eps$ of any representative of the germ $f$.
We shall make a similar abuse for vectors, matrices, polynomials whose
coefficients are germs.
We call \new{exponent} of $f\in\cont$:
\begin{equation}\label{puis-1}
\dex{f}\bydef \liminf_{\eps\to 0} \frac{\log | f(\eps)|}{\log \eps}
\in \R\cup\{+\infty\}\enspace .
\end{equation}

We have, for all $f,g\in \cont$ and $\lambda\in \C$,
\begin{align}
\dex{f+g}\geq & \min( \dex{f}, \dex{g})\enspace,\label{puis1}\\
\dex{fg}\geq & \dex{f}+ \dex{g}\enspace,\label{puis2}
\end{align}
with equality in~\eqref{puis1} if $\dex{f}\neq \dex{g}$ and
equality in~\eqref{puis2} if the liminf in the definition
of $\dex{f}$ or $\dex{g}$ is a limit.
Thus, $f\mapsto \dex{f}$ is ``almost'' a morphism $\cont\to\rmin$.
In the sequel exponents will be considered as elements of $\rmin$,
so that~\eqref{puis2} will be written as $\dex{fg}\geq \dex{f} \dex{g}$.
An element $f\in \cont$ is invertible  if, and only if,
$\dex{f}\neq \zero$ (or equivalently, if there exists
a positive constant such that $|f(\eps)|\geq \eps^k$).
If $f$ is invertible,  its inverse is the map 
$f^{-1}:\eps\mapsto f(\eps)^{-1}$ and we have $\dex{f^{-1}}\leq \dex{f}^{-1}$
with equality if, and only if, the liminf in the definition
of $\dex{f}$ is a limit.

We shall say that $f\in \cont$ has a \new{first order asymptotics} if
\begin{align}\label{puis4}
f(\eps)\sim a \eps^A, \qquad \mrm{when } \eps\to 0^+\enspace ,
\end{align}
with either $A\in \R$ and $a\in \C\setminus\{0 \}$,
or $A=+\infty$ and $a\in \C$.
In the first case,~\eqref{puis4} means that
$\lim_{\eps\to 0} \eps^{-A}f(\eps) =a$, in the second case,~\eqref{puis4} 
means that $f=0$ (in a neighbourhood of $0$).
We have:
\begin{align}
f(\eps)\sim a\eps^A & \implies  \dex{f}=A \enspace,\label{puis0}
\end{align}
and the liminf in~\eqref{puis-1} is a limit.
We shall also need an equivalence notion slightly weaker 
than $\sim$. If $f\in \cont$, $a\in \C$ and $A\in \rmin$, we write
\begin{equation}
f(\eps) \simeq a\eps^A
\label{e-f}
\end{equation}
if $f(\eps) = a\eps^A+o(\epsilon^A)$.
If $A\in\R$, this means that $\lim_{\eps\to 0} \eps^{-A}f(\eps) =a$.
If $A=+\infty$, this means by convention that $f=0$.
If $a\neq 0$ or $A=+\infty$, 
then $f(\eps) \simeq a\eps^A$ if, and only if, $f(\eps) \sim a\eps^A$
and in that case $\dex f =A$.
In general,
\begin{align}
f(\eps)\simeq a\eps^A & \implies  \dex{f}\geq A \enspace.\label{puis0-1}
\end{align}

Conversely, $\dex{f}>A\implies f(\eps)\simeq 0 \eps^A$. Of course, 
in~\eqref{e-f}, $a\eps^A$ must be viewed as a formal
expression, for the equivalence to be meaningful when $a=0$ and $A\in\R$.
In~\eqref{puis0}, however, $a\eps^A$ can  be viewed either as a formal
expression or as an element of $\cont$.

Throughout the paper, we consider a matrix $\sA\in \cont^{n\times n}$ 
and we shall assume that the entries $(\sA_\eps)_{ij}$ of $\sA_\eps$
have asymptotics of the form:
\begin{align}
(\sA_\eps)_{ij}\simeq a_{ij} \eps^{A_{ij}},\;\label{assump1} 
\mrm{for some matrix }a=(a_{ij})\in \C^{n\times n},\\
\mrm{and for some irreducible matrix }
A=(A_{ij})\in \rmin^{n\times n}.\nonumber
\end{align}
(The case where $A$ is reducible is a straightforward
extension.)
Under rather general circumstances (see Section~\ref{sec-newton}),
the eigenvalues $\sL_\eps^1,\ldots ,\sL_\eps^n$ of $\sA_\eps$ 
belong to $\cont$ and have first order asymptotics:
\begin{align}\label{asymp-l}
\sL_\eps^i\sim \lambda_i \eps^{\Lambda_i}\enspace .
\end{align}

We next relate the sequence $(\Lambda_1,\ldots,\Lambda_n)$
with two sequences constructed by using only the information
on the exponents of the entries $(\sA_\eps)_{ij}$ of the matrix $\sA_\eps$ 
given by the $A_{ij}$.

\subsection{First order Newton-Puiseux theorem and min-plus polynomials}\label{sec-newton}
The usual way to compute the $\Lambda_i$ in~\eqref{asymp-l}
is to use the classical Newton-Puiseux theorem.
We state here a general first order
version of this theorem
in a way which illuminates the role of min-plus algebra.

For any formal polynomial with coefficients in $\cont$,
$\sP (\eps, \iY)= \sum_{j=0}^{n} \sP_j(\eps) \iY^j\in \cont [\iY]$,
we define the min-plus \new{polynomial of exponents}:
\[ \dex{\sP}\bydef \bigoplus_{j=0}^n \dex{\sP_j} \iY^j\in \rmin[\iY]
\enspace .\]
The transformation of ordinary polynomials to min-plus 
(or ``tropical'') polynomial by the map $\dexo$ 
is instrumental in works on amoebas
(for instance, a very similar definition is given in~\cite{speyer}).

Recall that to $P=\dex{\sP}$ is associated the polynomial
function $\eval{P}$ and the
convex formal polynomial $\divex{P}$, as in Section~\ref{sec-minpol}.
For instance, to $\sP=\iY^3 +\eps^5 \iY^2 - \eps^6 \iY+ \eps^{13}$
corresponds the formal min-plus polynomial
$P=\dex{\sP}= \iY^3 \oplus 5\iY^2 \oplus 6\iY + 13$
represented in Figure~\ref{minpoly}.

\begin{theorem}[First order Newton-Puiseux theorem]\label{th3}
Let $\sP= \sum_{j=0}^{n} \sP_j(\eps) \iY^j\in \cont [\iY]$
such that $\sP_n=1$.
The following assertions are equivalent:
\begin{enumerate}
\item\label{th3-1} There exist $\sY_1,\ldots, \sY_n\in \cont$
such that $\sY_1(\eps),\ldots, \sY_n(\eps)$
are the roots of $\sP(\eps, y)=0$ counted with multiplicities, and
$\sY_1,\ldots, \sY_n$ have first order  asymptotics,
$\sY_j(\eps)\sim y_j \eps^{Y_j}$ with $Y_1\leq \cdots \leq Y_n$;
\item\label{th3-2} There exist $p=\sum_{j=0}^{n} p_j\iY^j\in \C [\iY]$
and $P=\bigoplus_{j=0}^{n} P_j\iY^j\in \rmin [\iY]$ satisfying
$\sP_j(\eps)\simeq p_j \eps^{P_j}$, $j=0,\ldots , n$,
with $p_n=1$, $P_n=\unit$, $p_0\neq 0$ or $P_0=\zero$, and 
$p_{n-i}\neq 0$ for all $i\in\{1,\ldots, n-1\}$ such that 
$c_i<c_{i+1}$, where $(c_1\leq \cdots \leq c_n)=\corn{P}$.
\end{enumerate}
When these assertions hold, we have 
$\dex{\sP}\geq P$, $\divex{\dex{\sP}}=\divex{P}$, and
$\corn{\dex{\sP}}=\corn{P}=(c_1\leq \cdots \leq c_n)=
(Y_1\leq \cdots \leq Y_n)$.
Moreover, if $c\in\R$ is a root of $P$ with multiplicity $k$ and
$c_i=\cdots =c_{i+k-1}=c$, then $y_i,\ldots, y_{i+k-1}$ are precisely
the non-zero roots of the polynomial 
\[
p^{(i)}=\sum_{\scriptstyle 0\leq j\leq n\atop\scriptstyle \eval{P}(c)=P_j c^j} p_j\iY^j\in \C[\iY]
\enspace,
\]
counted with multiplicities.
\end{theorem}

The classical  Newton-Puiseux theorem applies to the case where
$\cont$ is replaced by the field of (formal, or convergent)
Puiseux series (a Puiseux series is of the form
$ \sum_{k=K}^{\infty}  a_k x^{k/s}$ with $a_k \in \C$,
$K\in \Z$ and $s\in \N\setminus\{0\}$),
and shows~\ref{th3-2}$\implies$\ref{th3-1} only.
In the classical statement of the theorem, 
the leading exponents $Y_i$, are, up to an inversion 
and change of sign, the slopes of the Newton polygon, and the 
polynomials $p^{(i)}$ are defined in terms of the 
edges of the polygon. 
Since, when $P=\dex{\sP}$, the graph of $\vex P$ is the 
symmetric, with respect to the main diagonal, of the Newton polygon, 
it follows from  Proposition~\ref{prop-prenewton} that the 
$Y_i$ and $y_i$ in Theorem~\ref{th3} coincide with the ones
that are defined classically.

Theorem~\ref{th3} is a ``precise
large deviation'' version of the Newton-Puiseux theorem:
we assume only the existence of asymptotic equivalents
for the coefficients of $\sP(\epsilon,\cdot)$,
and derive the existence of asymptotic equivalents
for the branches of $\sP(\epsilon,\cdot)$.
The Newton-Puiseux algorithm is sometimes presented
for asymptotic expansions, as in~\cite{dieu}.
However, the equivalence between the two assertions
of Theorem~\ref{th3} does not seem to be classical.
In particular, the asymptotics of some coefficients
may be only known as being negligible: we require that $p_i\neq 0$
only for those $i$ such that $(i,P_i)$ is an exposed
point of the epigraph of $P$.

\begin{proof}
We first prove~\ref{th3-1}$\implies$\ref{th3-2}.
Let $Q=(\iY\oplus Y_1)\cdots (\iY\oplus Y_n)$. Then, $Q=\divex{Q}$,
$\corn{Q}=(Y_1\leq \cdots \leq Y_n)$ and
$Q_{n-i}= Y_1\cdots Y_i$ for all $i=1,\ldots , n$.
Since $\sY_1(\eps),\ldots, \sY_n(\eps)$
are the roots of $\sP(\eps, y)=0$ counted with multiplicities, 
and $\sP_n=1$, it follows that
$\sP(\eps, \iY) =\prod_{i=1}^n (\iY- \sY_i(\eps)) $.
Hence, $(-1)^i\sP_{n-i}$ is the sum of all
products $\sY_{j_1}\cdots \sY_{j_i}$, where $j_1,\ldots, j_i$
are pairwise distinct elements of $\{1,\ldots , n\}$.
By the properties of ``$\simeq$'' (stability by addition and
multiplication), and since $\bigoplus_{j_1,\ldots , j_i}
 Y_{j_1}\cdots Y_{j_i}=Y_1\cdots Y_i
=Q_{n-i}$, we obtain that there exist $p_0,\ldots, p_{n-1}\in \C$
such that $\sP_{j}\simeq p_j \eps^{Q_j}$ for all $j=0,\ldots, n-1$.
Putting $p_n=1$, we also get $\sP_{n}=1\simeq p_n \eps^{Q_n}$ 
since $Q_n=\unit$.
When $i=1,\ldots , n-1$ is such that
$Y_i<Y_{i+1}$, $\sY_{1}\cdots \sY_{i}$ is the only leading term 
in the sum of all $\sY_{j_1}\cdots \sY_{j_i}$,
and then $p_{n-i}=(-1)^i y_1\cdots y_i\neq 0$.
Moreover, for $i=n$, either $Y_n\neq\zero$, which implies that
$p_{0}=(-1)^n y_1\cdots y_n\neq 0$, or $Y_n=\zero$, which implies
that $\sY_n=0$, $\sP_0=0$ and $Q_0=\zero$.
This shows that $(c_1,\ldots,c_n)=(Y_1,\ldots,Y_n)$
and $P=Q$ are as in Point~\ref{th3-2}.

The remaining part
of the theorem is obtained by a simple adaptation of the proof
of the classical Newton-Puiseux theorem.
When the $\sP_j$ are only assumed to be continuous functions
satisfying Point~\ref{th3-2} of the theorem, it
follows from~{\rm (\ref{puis1},\ref{puis2},\ref{puis0-1})}, 
that $\dex{\sP}\geq P$, and since $P\geq \divex{P}=(\iY\oplus c_1)\cdots
(\iY\oplus c_n)$, we get that $\dex{\sP}\geq (\iY\oplus c_1)\cdots
(\iY\oplus c_n)$.  In addition,
 from~\eqref{puis0} and Point~\ref{th3-2} of the theorem, we get that
$\dex{\sP}_{n-i}= P_{n-i}=\divex{P}_n c_1\cdots \cdots c_i$
for all $i\in\{0,n\}\cup\set{i\in \{1,\ldots , n-1\}}{c_i< c_{i+1}}$, 
hence  Lemma~\ref{carac-ci} yields $\divex{\dex{\sP}}=\divex{P}$,
therefore, $\corn{\dex{\sP}}=\corn{P}=(c_1\leq\cdots \leq c_n)$.
Moreover, the first step of the Puiseux algorithm shows
that, for all roots $c\neq\zero$ of $P$ with multiplicity $k$,
there are exactly $k$ continuous branches with leading exponent $c$.
Indeed, when $c=c_i=\cdots =c_{i+k-1}\neq \zero$,
the change of variable $y=z \eps^{c}$, and the division of $\sP$
by $\eps^{\eval{P}(c)}$, transforms the equation $\sP(\eps,y)=0$ into
an equation $\sQ(\eps,z)=0$, where $\sQ(\cdot,z)$ extends 
continuously to $0$ with $\sQ(0,z)=p^{(i)}(z)$.
Since $\eval{P}(c)=P_j c^j$ implies that $n-i-k+1\leq j\leq n-i+1$,
and since  either $i-1=0$ or $c_{i-1}<c_i$,
we get that $p_{n-i+1}\neq 0$, hence $\deg p^{(i)}=n-i+1$.
Similarly, we have either $i+k-1=n$ or $c_{i+k-1}<c_{i+k}$.
In the second case, we get $p_{n-i-k+1}\neq 0$, thus $\val p^{(i)}=n-i-k+1$.
In the first case, $i+k-1=n$, $c=c_n$,  and $p_{0}\neq  0$ or $P_{0}=\zero$.
Since $P_{0}=\zero$ implies $c=c_n=\zero$, which contradicts our
assumption, we must have $p_0\neq 0$, hence again $\val p^{(i)}=n-i-k+1$.
Hence, $\deg p^{(i)}-\val p^{(i)}=k$ and 
the conclusion is obtained by the standard Lemma~\ref{cont-roots}
below.
Finally, if $c=\zero$ is a root with multiplicity $k$, then $\val P=k$,
$c_{n-k}<c_{n-k+1}=\zero$, and $p_{k}\neq 0$. This implies that 
(for all $\eps>0$ in a
neighbourhood of $0$) $\sP(\eps,\cdot)$ is a polynomial with valuation 
$k$, hence it has $0$ as a root with multiplicity $k$.
\end{proof}
\begin{lemma}\label{cont-roots}
Let $\sQ(\eps,\iY)=\sum_{i=0}^n \sQ_j(\eps) \iY^j$, where the $\sQ_j$ are 
continuous functions of $\eps\in [0,\eps_0)$ and let $m=\deg \sQ(0,\cdot)$.
Then, for any open ball $B$ containing the roots of $\sQ(0,\cdot)$,
there are $m$ continuous branches $\sZ_1,\ldots ,\sZ_m$ defined in 
some interval $[0,\eps_1)$, with $0<\eps_1\leq\eps_0$, such that  
$\sZ_1(\eps),\ldots ,\sZ_m(\eps)$ are exactly the roots of 
$\sQ(\eps,\cdot)$ in $B$ counted with multiplicities.
Moreover, the roots of $\sQ(\eps,\cdot)$ that are outside $B$
tend to infinity when $\eps$ goes to $0$.
\end{lemma}
\begin{proof}
We only sketch the proof, which is classical.
By the Cauchy index theorem, 
if $\gamma$ is any circle in $\C$
containing no roots of $\sQ(\epsilon,\cdot)$,
the number of roots of $\sQ(\epsilon,\cdot)$
inside $\gamma$ 
is $(2\pi i)^{-1}\int_\gamma \partial_z \sQ(\epsilon,z)(\sQ(\epsilon,z))^{-1}\,\mrm{d}z$.
By continuity of $\eps\mapsto \sQ(\eps,\cdot)$, 
the number of roots of $\sQ(\eps',\cdot)$
inside $\gamma$ (counted with multiplicities) is constant for $\eps'$ in
some neighbourhood of $\eps$. Taking $B$ as in the lemma,
$\gamma=\partial B$, and $\eps=0$, we get exactly $m$  roots of $\sQ(\eps',\cdot)$
in $B$ for $\eps'$ in some interval $[0,\eps_1)$.
Consider now a ball $B_R\supset B$ of radius $R$. 
For $\eps'$ small enough, the number of roots of $\sQ(\eps',\cdot)$ in
either $B_R$ or $B$ is equal to $m$, hence
any root of $\sQ(\eps',\cdot)$ outside $B$ must be outside $B_R$.
This shows that the roots of $\sQ(\eps',\cdot)$ 
that do not belong to $B$ go to infinity, when $\eps'\to 0$.
Finally, by taking small balls around each root of $\sQ(\eps,\cdot)$,
with $0\leq \eps< \eps_1$, we see that the map which sends $\eps$
to the unordered $m$-tuple of roots of $\sQ(\eps,\cdot)$ that belong to $B$,
is continuous on $[0,\eps_1)$.
By a selection theorem for unordered $m$-tuples depending continuously on a
real parameter (see for instance~\cite[Ch. II, Section~ 5, 2]{kato}),
we derive the existence of the $m$ continuous branches $\sZ_1,\ldots ,\sZ_m$.
\end{proof}

Theorem~\ref{th3} says that ``the leading exponents of the roots are
the min-plus roots''. 
\begin{example}
Consider $\sP(\eps,\iY)=\iY^3 +\eps^5 \iY^2 -\eps^6 \iY + \eps^{13}$. The 
min-plus polynomial $P=\dex{\sP}$ is the one of Figure~\ref{minpoly},
hence its roots are $c_1=c_2=3$ and $c_3=7$.
We have $p^{(1)}=p^{(2)}=\iY^3-\iY$ and $p^{(3)}=-\iY+1$.
Hence, $\sP$ has $3$ continuous branches around $0$
with first order asymptotics: $\sY_1\sim \eps^3$, $\sY^2\sim-\eps^3$
and $\sY_3\sim \eps^7$. 
Theorem~\ref{th3} states in particular that we need
not know the asymptotic expansions of all
the coefficients of $\sP(\eps,\iY)$:
for instance,
if $\sP(\eps,\iY)=\iY^3 +o(\eps^3) \iY^2 -\eps^6 \iY + \eps^{13}$,
the polynomials $P$ and $p^{(1)},p^{(2)},p^{(3)}$ are unchanged, so that
we still have $3$ continuous branches with the same
asymptotics has above. 
\end{example}
\begin{remark}\label{rem-puis}
If $\sA\in\cont^{n\times n}$ satisfies~\eqref{assump1}, the characteristic 
polynomial of $\sA_\eps$, $\sP(\eps,\iY)= \det (\iY I- \sA_\eps)$ 
is an element of $\cont[\iY]$, since $\cont$ is a ring. 
Applying Theorem~\ref{th3} to $\sP$, we can obtain, under
some additional assumptions,
first order asymptotics for the eigenvalues of $\sA_\eps$.
The difficulty is that the coefficients $\sP_j$ of $\sP$ need
not have first order asymptotics (even if $a_{ij}\neq 0$ for all
$i,j$) due to cancellations. 
Of course if the coefficients 
of $\sA_\eps$ have Puiseux series expansions in $\eps$, the
$\sP_j$ also  have Puiseux series expansions in $\eps$
and a fortiori first order asymptotics.
However, if we only assume that $\sA\in\cont^{n\times n}$ 
satisfies~\eqref{assump1}, we obtain that
 the $\sP_j$ satisfy the conditions $\sP_n= 1$
and $\sP_j(\eps)\simeq p_j \eps^{P_j}$ for some exponents 
$P_j\in\rmin$ computed using 
the exponents $A_{ij}$ (see Section~\ref{sec-min-plus-charac}).
Hence,  if the eigenvalues of $\sA_\eps$ have first order asymptotics,
Theorem~\ref{th3} gives the exponents of these asymptotics as
a function of the $P_j$. 
\end{remark}

\subsection{Majorisation inequalities for roots of min-plus polynomials}
\label{sec-min-plus-charac}
The permanent of a matrix with coefficients in an arbitrary
semiring $(S,\oplus,\otimes)$ can be defined as usual:
\begin{align*}
\perm (A)&= \bigoplus_{\sigma\in \Sym_n}
\bigotimes_{i=1}^n A_{i\sigma (i)} \enspace,
\end{align*}
where $\Sym_n$ is the set of permutations of $\{1,\ldots, n\}$.
In particular, for any matrix $A\in \rmin^{n\times n}$, 
\begin{align*}
\perm (A) =\min_{\sigma\in \Sym_n}
\sum_{i=1}^n A_{i\sigma (i)} \enspace ,
\end{align*}
and the \new{formal characteristic polynomial} of $A$
is the polynomial 
\[
\perm (\iY I \oplus A)= \bigoplus_{\sigma\in \Sym_n}
\bigotimes_{i=1}^n (\iY \delta_{i\sigma (i)} \oplus A_{i\sigma (i)})
\in \rmin[\iY] \enspace,
\]
where $I$ is the identity matrix,
and $\delta_{ij}=\unit$ if $i=j$ and $\delta_{ij}=\zero$ otherwise.
The associated min-plus polynomial function will be called
the \new{characteristic polynomial function} of $A$.

We next assume that $\sA\in \cont^{n\times n}$ satisfies~\eqref{assump1}
and that the eigenvalues $\sL^i_\eps$ ($i=1,\ldots , n$) 
of $\sA_\eps$ have first order asymptotics,
 $\sL^i_\eps\sim \lambda_i \eps^{\Lambda_i}$. We relate,
in that case, the $\Lambda_i$ with the roots
of the characteristic polynomial of $A$.

We need first to recall the classical definition of weak majorisation
(see~\cite{MAR} for more background).
\begin{definition}
Let $u,v\in \rmin^n$. Let $u_{(1)}\leq \cdots \leq u_{(n)}$
(resp.\ $v_{(1)}\leq \cdots \leq v_{(n)}$) denote the components
of $u$ (resp.\ $v$) in increasing order.
We say that $u$ is \new{weakly (super) majorised} by $v$,
and we write $u\weakm v$, if the following conditions hold:
\[  u_{(1)}\cdots  u_{(k)}
\geq  v_{(1)} \cdots   v_{(k)}  \quad \forall k=1, \ldots, n \enspace . \]
\end{definition}

In fact, the weak majorisation relation is only defined
in~\cite{MAR} for vectors of $\R^n$. Here, it is convenient
to define this notion for vectors of $\rmin^n$.
We also used the min-plus notations for homogeneity with the rest 
of the paper.
The following lemma states a useful
monotonicity property of the map which associates to a
formal min-plus polynomial $P$ its sequence
of roots, $\corn{P}$. 
\begin{lemma}\label{minpoly-maj}
Let $P, Q\in \rmin[X]$ be two formal polynomial of degree $n$.
Then,
\begin{equation}
P\geq Q \mrm{ and } P_n=Q_n \implies \corn{P}\weakm \corn{Q}
\enspace .
\end{equation}
\end{lemma}
\begin{proof}
From $P\geq Q$, we deduce $\divex{P}\geq \divex{Q}$.
Let $\corn{P}= (c_1(P)\leq\cdots \leq c_n(P))$
and 
$\corn{Q}=(c_1(Q)\leq\cdots\leq c_n(Q))$ denote the sequence of roots
of $P$ and $Q$, respectively.
Using $\divex{P}\geq \divex{Q}$, $\divex{P}_n=P_n=Q_n=\divex{Q}_n$ and 
\eqref{corners0}, we get
$c_1(P)\cdots c_k(P)=\divex{P}_{n-k}(\divex{P}_{n})^{-1}
\geq \divex{Q}_{n-k}(\divex{Q}_{n})^{-1}=c_1(Q)\cdots c_k(Q)$,
for all $k=1,\ldots , n$, that is $\corn{P}\prec^{\rm w} \corn{Q}$.
\end{proof}

We shall also need the following notion of genericity.
We will say that a property $\ptyP(y)$ depending on the variable $y=(y_1,\ldots
, y_n)\in \C^n$ holds for \new{generic values} of $y$ if the set 
of elements  $y\in \C^{n}$ such that the property $\ptyP(y)$
is false is a proper algebraic variety.
This means that there exists $Q\in\C[\iY_1,\ldots, \iY_n]\setminus\{0\}$ 
such that $\ptyP(y)$ is false if $Q(y)=0$.
When the parameter $y$ will be obvious, we shall simply say that
$\ptyP$ is \new{generic} or holds \new{generically}.
It is clear that if $\ptyP_1$ and $\ptyP_2$ are both generic, then
``$\ptyP_1$ and $\ptyP_2$'' is also generic.

Since any polynomial $q=\sum_{i_1,\ldots, i_n\in\N} 
q_{i_1,\ldots, i_n} \iY_1^{i_1}\cdots \iY_n^{i_n} 
\in \C[\iY_1,\ldots, \iY_n]$ in $n$ indeterminates
can be seen as
an element of  $\cont [\iY_1,\ldots, \iY_n]$ whose coefficients are constant 
with respect to $\eps$, we have:
\begin{align}\label{qQ}
\dex{q}=\bigoplus_{\scriptstyle i_1,\ldots, i_n\in\N\atop\scriptstyle q_{i_1,\ldots, i_n}\neq 0}
\iY_1^{i_1}\cdots \iY_n^{i_n} \enspace \in \rmin[\iY_1,\ldots, \iY_n]
\enspace .\end{align}
We also define, for all $Y\in\rmin^n$:
\begin{align}\label{qsat}
q^{\sat}_Y(\iY):=\sum_{\scriptstyle i_1,\ldots, i_n\in\N\atop\scriptstyle \dex{q}(Y_1,\ldots , Y_n)=
Y_1^{i_1}\cdots Y_n^{i_n}} q_{i_1,\ldots, i_n} \iY_1^{i_1}\cdots
\iY_n^{i_n} \in \C[\iY_1,\ldots, \iY_n] \enspace .
\end{align}

The following result is clear from the above definitions of $\dex{q}$
and $q^\sat_Y$, since when $y\neq 0$, $\sY\simeq y \eps^Y\iff
\sY\sim y \eps^Y$.
\begin{lemma}\label{gener}
Let $q\in \C[\iY_1,\ldots, \iY_n]$ and  let $Q=\dex{q}$ and
$q^\sat_Y$ be defined by~\eqref{qQ} and~\eqref{qsat}, respectively.
Let $\sY\in \cont^n$, $y\in \C^n$ and $Y\in \rmin^n$ be 
such that $\sY_i\simeq y_i \eps^{Y_i}$ for $i=1,\ldots , n$. Then, 
\begin{align}\label{simpoly}
 q(\sY_1,\ldots, \sY_n)&\simeq q^{\sat}_Y(y) \eps^{Q(Y)}\enspace ,
\end{align}
and for any fixed $Y$, we have an equivalence $\sim$ in~\eqref{simpoly}
for generic values of $y\in \C^n$.
\end{lemma}

\begin{theorem}\label{th-gen=}
Let  $\sA\in \cont^{n\times n}$ satisfy~\eqref{assump1}.
Assume that the eigenvalues $\sL^1_\eps,\ldots , \sL^n_\eps$
of $\sA_\eps$ (counted with multiplicities) have first order asymptotics,
 $\sL^i_\eps\sim \lambda_i \eps^{\Lambda_i}$,
and denote by $\Lambda=(\Lambda_1\leq \cdots \leq \Lambda_n)$ 
the sequence of their exponents (counted with multiplicities). Let 
$\Gamma=(\gamma_1\leq\cdots \leq \gamma_n)$ be the sequence of roots 
of the min-plus characteristic polynomial of $A$.
Then,
\begin{align}\label{lwg}
\Lambda  \weakm\Gamma \enspace ,
\end{align}
and for generic values of $a=(a_{ij})\in \C^{n\times n}$, 
$\Lambda=\Gamma$.
\end{theorem}
\begin{proof}
Since $\sA=\sA_\epsilon\in \cont^{n\times n}$,
the characteristic polynomial of $\sA$,
$\sQ(\eps,\iY):=\det (\iY I- \sA_\eps)$ belongs to $\cont[\iY]$.
Let $Q=\dex{\sQ}\in \rmin[\iY]$ and let
$P=\perm (\iY I\oplus A)\in  \rmin[\iY]$
be the min-plus characteristic polynomial of $A$.
By~{\rm (\ref{puis1},\ref{puis2})}, 
$\dex{\sQ}\geq \perm(\iY I\oplus \dex{-\sA})$ and by~\eqref{assump1}
and~\eqref{puis0-1},
$\dex{-\sA}\geq A$. It follows that $Q=\dex{\sQ}\geq 
\perm(\iY I\oplus A)=P$. Hence, from Lemma~\ref{minpoly-maj}, 
we get that $\corn{Q}\weakm \corn{P}=\Gamma$.
Moreover, by Theorem~\ref{th3} applied to $\sQ$,
we get that  $\corn{Q}=\Lambda$, which finishes the proof of~\eqref{lwg}.

Let us show the genericity of the equality $\Lambda=\Gamma$.
For all $a\in \C^{n\times n}$, we consider the $k$-th trace of $a$:
\begin{align*}
\tr_{k}(a)=\sum_{J\subset\{1,\ldots , n\},\, \card{J}= k}\left(
\sum_{\sigma\in \Sym_J} \sgn (\sigma) \prod_{j\in J}  a_{j\sigma (j)} 
\right) \enspace .\end{align*}
For all $A\in \rmin^{n\times n}$, we also set
\begin{equation}\label{trmin}
\tr_k(A)=
\bigoplus_{J\subset\{1,\ldots , n\},\, \card{J}= k}\left(
\bigoplus_{\sigma\in \Sym_J} \bigotimes_{j\in J}  A_{j\sigma (j)} \right)
\enspace .
\end{equation}
Then, the coefficients of $\sQ$ are given by 
 $\sQ_k(\eps)=(-1)^{k} \tr_{n-k} (\sA_\eps)$, for $k=0,\ldots , n-1$
and $\sQ_n=1$.
The coefficients of $P$ are given by  $P_k=\tr_{n-k} (A)$, for $k=0,\ldots ,
 n-1$ and $P_n=\unit$.
By Lemma~\ref{gener}, we obtain that for any fixed (irreducible)
matrix $A\in\rmin^{n\times n}$, and any $\sA\in\cont^{n\times n}$
satisfying~\eqref{assump1} with $a\in\C^{n\times n}$ and $A$,
$\tr_k(\sA_\eps)\sim(\tr_k)^\sat_A(a) \eps^{\tr_{k}(A)}$
for generic values of $a\in \C^{n\times n}$.
In particular, generically, 
$\sQ_k(\eps)$ has first order asymptotics and
$\dex{\sQ_k}=P_k$, for all $k=0,\ldots , n$. This implies that
$Q=P$, thus $\Lambda=\corn{Q}=\corn{P}=\Gamma$, generically.
\end{proof}
\begin{remark}
Since a result of Burkard and Butkovi\v{c}~\cite{bb02} shows
that we can compute the min-plus characteristic
polynomial function of a matrix in polynomial time (by solving $O(n)$
assignment problems), Theorem~\ref{th-gen=} shows that
the sequence $\Lambda$ of generic exponents of the eigenvalues
can be computed in polynomial time. See also~\cite{murota,abg04}. 
\end{remark}
\section{Critical values of min-plus matrices}
\subsection{Schur complements and generalised circuit means} \label{sec-graph} 
We now construct another sequence $\beta=(\beta_1\leq \cdots \leq \beta_n)$
using eigenvalues of min-plus matrices.
First, we build by induction a finite
sequence of min-plus square matrices
$A_{\ell}$ and scalars $\alpha_{\ell}\in \R$, 
for $1\leq {\ell}\leq k$,
together with a partition
$C_1\cup \cdots \cup C_k=\{1,\ldots,n\}$.

We start with $A_1=A$. Then, for all ${\ell}\geq 1$, we define
\begin{equation}
\alpha_{\ell}=\rhomin(A_{\ell})
\label{e-alpha}
\end{equation}
and we take for $C_{\ell}$ the set of critical
nodes of $A_{\ell}$. We build, as long as
$C_1\cup \cdots \cup C_{\ell}\neq \{1,\ldots,n\}$,
the min-plus Schur complement:
\[
A_{{\ell}+1}= \sch(C_{\ell}, \alpha_{\ell},A_{\ell}) \enspace .
\]
Due to the irreducibility of $A$, Lemma~\ref{lem-schur} shows that
$A_{\ell}$ is irreducible, so that $C_{\ell}\neq \emptyset$.
Hence, the algorithm stops at some index $k\leq n$.
By Proposition~\ref{prop1}, we get that $\alpha_1 < \cdots < \alpha_k$.
We call  $\alpha_1,\ldots ,\alpha_k$ the \new{critical values} of $A$.
We define the \new{multiplicity} of the critical value $\alpha_{\ell}$ as
$\card C_{\ell}$. Repeating each  critical value with its multiplicity, we obtain
a sequence $\beta=(\beta_1\leq \cdots \leq \beta_n)$ which will be called
the \new{sequence of critical values counted with multiplicities}.

Let us give now a graph interpretation of the exponents $\alpha_{\ell}$.
We set $C^0=\emptyset$ and, for all ${\ell}=1,\ldots , k$, 
\[
C^{\ell}=C_1\cup\ldots\cup C_{{\ell}},
\qquad 
N^{\ell}=\{1,\ldots,n\}\setminus C^{{\ell}-1}
\enspace.
\]
For all paths $p$ of $G(A)$ and all ${\ell}=1,\ldots , k$, we use the 
notations of Section~\ref{sec-schur} and:
\begin{align*}
|p|^{\ell}_A :=& |p|_A-\alpha_1|p|_{C_1} -\cdots -\alpha_{{\ell}-1}|p|_{C_{{\ell}-1}} 
\enspace , \\
|p|^{\ell}:=& |p|-|p|_{C_1}-\cdots-|p|_{C_{{\ell}-1}}=|p|_{N^{\ell}}\enspace .
\end{align*}
\begin{proposition}\label{lem-alpha}
The numbers $\alpha_{\ell}$ defined in \eqref{e-alpha} satisfy:
\begin{equation}
\alpha_{\ell}=\min \frac{|c|^{\ell}_A}{|c|^{\ell}}\enspace,
\label{e-r2}
\end{equation}
where the minimum is taken over all circuits $c$ in $G(A)$
which are not included in $C^{{\ell}-1}$.
Moreover, $c$ is a critical circuit of $A_{\ell}$ if, and only if, there exists
a circuit $c'$ of $G(A)$ such that $c'\cap N^{\ell}=c$ and $c'$
minimises~\eqref{e-r2}.
\end{proposition}
\begin{proof}
Using repetitively Lemma~\ref{lem-schur}, we get that for
all circuits $c$ of $G(A_{\ell})$, $|c|_{A_{\ell}}=\min |c'|^{\ell}_A$,
where the minimum is taken over all circuits $c'$ of $G(A)$ such that
$c'\cap N^{\ell}=c$.
By the same arguments as in the proof of Proposition~\ref{prop1},
we deduce the assertions of Proposition~\ref{lem-alpha}.
\end{proof}

Note that, as for Proposition~\ref{prop1},
 if $c$ is included in $C^{{\ell}-1}$,
that is if the denominator in~\eqref{e-r2} is zero, 
the numerator is necessarily nonnegative
(by definition of $\alpha_{{\ell}-1}$).

We say that a circuit $c$ of $G(A)$ is 
a \new{critical circuit of order $\ell$}
if $|c|^{\ell}_A=\alpha_{\ell} |c|^{\ell}$.
We call \new{critical graph of order
${\ell}$} the graph $G^c_{\ell}(A)$ whose nodes and
arcs belong to critical circuits of order ${\ell}$. 
Of course, $G^c(A)=G^c_1(A)$.

\begin{proposition}
We have 
\begin{equation}
G^c_{\ell}(A)\subset G^c_{{\ell}+1}(A)\quad {\ell}=1,\ldots , k-1 \enspace .
\label{e-nested}
\end{equation}
which means that the nodes and arcs of $G^c_{\ell}(A)$
belong to $G^c_{{\ell}+1}(A)$.
\end{proposition}
\begin{proof}
If $c$ is a critical circuit of order ${\ell}$, then by definition
$|c|^{\ell}_A=\alpha_{\ell} |c|^{\ell}$.
If in addition $|c|^{\ell}=0$, then $c\cap N^{\ell}=\emptyset$, whence $|c|_{C_{\ell}}=0$ and
$|c|^{{\ell}+1}=0$. It follows that $|c|^{{\ell}+1}_A=|c|^{\ell}_A-\alpha_{\ell} |c|_{C_{\ell}}=
0= \alpha_{{\ell}+1} |c|^{{\ell}+1}$, thus $c$ is a critical circuit of order ${\ell}+1$.
Otherwise, if $|c|^{\ell}\neq 0$, then $c$ minimises~\eqref{e-r2} and 
since, by the arguments of the proof of Proposition~\ref{prop1},
$|c\cap N^{\ell}|_{A_{\ell}}\leq |c|^{\ell}_A$, we obtain that $c'=c\cap N^{\ell}$ is
a critical circuit of $A_{\ell}$. By definition of $C_{\ell}$, we get that
the nodes of $c'$ belong to $C_{\ell}$, thus the nodes of $c$ belong
to $C^{\ell}$, which shows  $|c|^{{\ell}}=|c|_{C_{\ell}}$ or $|c|^{{\ell}+1}=0$.
Since $|c|^{\ell}_A=\alpha_{\ell} |c|^{\ell}$, we get $|c|^{{\ell}+1}_A=0=|c|^{{\ell}+1}$,
and $c$ is a critical circuit of order ${\ell}+1$.
\end{proof}

Let $\ell\in\{1,\ldots , k\}$ and let $D_{\ell}$ denote the min-plus diagonal
matrix such that $(D_{\ell})_{jj} =\alpha_m$ if $j\in C_m$ with $m< {\ell}$,
and $(D_{\ell})_{jj} =\alpha_{\ell}$ if $j\in N^\ell$.
For instance, if $n=3$, $C_1=\{1\}$, $C_2=\{2,3\}$, 
$\alpha_1=2$ and $\alpha_2=4$, then $D_1=\diag(2,2,2)$ and
$D_2=\diag(2,4,4)$. We set 
\[ \hat{A}_{\ell} = D_{\ell}^{-1} A 
\enspace .
\]
We also set 
\[G_{\infty}^c(A)=G^c_k(A), \quad
\ainf =\hat{A}_k,\quad \mrm{ and }\dinf=D_k
\enspace .
\]
\begin{lemma}\label{aiachapi}
We have $A_{\ell}=\alpha_{\ell} \sch(C^{\ell-1},\hat{A}_{\ell})$, for ${\ell}=1,\ldots, k$.
\end{lemma}
\begin{proof}
We prove the lemma by induction on ${\ell}=1,\ldots , k$.
Since $\hat{A}_1=\alpha_1^{-1} A$ and $A_1=A$, we get $A_1=\alpha_1
\hat{A}_1$.
If $A_{\ell}=\alpha_{\ell} \sch(C^{\ell-1},\hat{A}_{\ell})$, then using~\eqref{schur-homo} 
and~\eqref{e-ident}, we get 
\begin{align*}
A_{{\ell}+1}&=\sch(C_{\ell},\alpha_{\ell}, A_{\ell})=\alpha_{\ell} \sch(C_{\ell}, \alpha_{\ell}^{-1} A_{\ell})\\
&= \alpha_{\ell} \sch(C_{\ell}, \sch(C^{\ell-1},\hat{A}_{\ell}))=\alpha_{\ell} \sch(C^{\ell}, \hat{A}_{\ell})\\
& = \alpha_{\ell} \sch(C^{\ell}, D_{\ell}^{-1} D_{{\ell}+1} \hat{A}_{{\ell}+1})\enspace .
\end{align*} 
Since $(D_{\ell}^{-1}D_{\ell+1})_{jj}=\unit$ for $j\in C^\ell$,
and $(D_{\ell}^{-1}D_{\ell+1})_{jj}=\alpha_\ell^{-1}\alpha_{\ell+1}$ 
otherwise, it follows from~\eqref{e-def-mschur} that
$A_{\ell+1}=\alpha_{\ell+1}\sch(C^\ell, \hat{A}_{\ell+1})$.
\end{proof}
\begin{proposition}\label{lem-3}
For all $1\leq {\ell}\leq k$, we have
$G^c_{\ell}(A)=G^c(\hat{A}_{\ell})$, $\hat{A}_{\ell}$ has
min-plus eigenvalue $\unit$, and
the set of critical nodes of $\hat{A}_{\ell}$ is $C^{\ell}$.
Moreover, $G_{\ell}^c(A)$  and $\restr{G^c(\ainf)}{C^{\ell}}$
(that is the restriction of $G^c(\ainf)$ to the nodes of $C^{\ell}$)
have the same strongly connected components.
In particular, $G^c_\infty (A)=G^c(\ainf)$
and all the nodes of $\{1,\ldots, n\}$ are critical for $\ainf$.
\end{proposition}

\begin{proof}
For all circuits $c$ and for all ${\ell}=1,\ldots , k$, we get by 
Proposition~\ref{lem-alpha},
$|c|^{\ell}_A\geq \alpha_{\ell} |c|^{\ell}$. Since, for all circuits
$|c|_{\hat{A}_{\ell}}= |c|_A-\alpha_1 |c|_{C_1}-\cdots -\alpha_{{\ell}-1} |c|_{C_{{\ell}-1}}
-\alpha_{\ell} |c|_{N^{\ell}}=|c|^{\ell}_A-\alpha_{\ell} |c|^{\ell}$, we get that
$\rho(\hat{A}_{\ell})\geq 0$.
Moreover, $c$ is a critical circuit of order ${\ell}$ if, and only if,
$|c|^{\ell}_A=\alpha_{\ell} |c|^{\ell}$, which is equivalent to $|c|_{\hat{A}_{\ell}}= 0$.
This shows that $\rho(\hat{A}_{\ell})= 0$ and that $c$ is a 
critical circuit of order ${\ell}$ if, and only if, $c$ is a critical circuit
of $\hat{A}_{\ell}$. It follows that $G_{\ell}^c(A)=G^c(\hat{A}_{\ell})$.
Since by Proposition~\ref{lem-alpha}, any critical circuit of $A_\ell$
is of the 
form $c'\cap N^{\ell}$ where $c'$ is a critical circuit of order ${\ell}$, the set 
$C_{\ell}$ of nodes of $G^c(A_{\ell})$ is included in the set of nodes of
$G^c_{\ell}(A)$. Using~\eqref{e-nested}, we get by induction
that $C^{\ell}$ is included in the set of nodes of $G^c_{\ell}(A)$. 
Conversely, since any critical circuit $c'$ of order ${\ell}$ is such that
$c'\cap N^{\ell}$ is a critical circuit of $A_{\ell}$, and since the set of critical 
nodes of $A_{\ell}$ is $C_{\ell}$, the set of nodes of $G^c_{\ell}(A)$ is included
in $(\{1,\ldots , n\}\setminus N^{\ell})\cup C_{\ell}=C^{\ell}$, hence is equal to $C^{\ell}$.
Finally it is clear that, by definition of $\hat{A}_{\ell}$,
$G^c(\hat{A}_{\ell})\subset G^c(\ainf)$, and since
its set of nodes is $C^{\ell}$, we get
$G^c(\hat{A}_{\ell})\subset \restr{G^c(\ainf)}{C^{\ell}}$ .
Conversely, since the restrictions of $\ainf$ and $\hat{A}_{\ell}$
to $C^{\ell}\cart C^{\ell}$ are equal and since $\rho(\hat{A}_{\ell})=\
\rho(\ainf)=0$, any critical circuit of $\ainf$ with nodes in $C^{\ell}$ 
is critical for $\hat{A}_{\ell}$.
It follows that the strongly connected components of $G^c(\hat{A}_{\ell})$
and $G^c(\ainf)$ are equal.
\end{proof}

\begin{example}\label{ex-critic}
To illustrate the computation of the critical values,
consider
\begin{align}
A=
\left[\begin{array}{ccccc}
\infty & 0 & \infty & \infty\\
0 & \infty & 1  & \infty\\
1 & \infty &\infty & 2 \\
\infty & \infty &4  & 5 
\end{array}\right]
\enspace .
\label{e-ex-critic}
\end{align}
We have $\alpha_1=0$, and the critical graph of $A$ is
composed of the circuit $(1\to 2\to 1)$. Thus, $C_1=\{1,2\}$.
We have
\begin{align*}
A_2&=\sch(C_1,\alpha_1, A)\\
&= 
\left[\begin{array}{cc}
\infty & 2\\
4  &  5
\end{array}
\right]
\oplus 
\left[\begin{array}{cc}
1 & \infty\\
\infty & \infty\\
\end{array}
\right]
\left[\begin{array}{cc}
\infty & 0\\
0 & \infty
\end{array}\right]^*
\left[\begin{array}{cc}
\infty & \infty \\
1  & \infty 
\end{array}\right]
= 
\left[\begin{array}{cc}
2 & 2\\
4  &  5
\end{array}\right]\enspace .
\end{align*}
Hence, $\alpha_2=\rhomin(A_2)=2$, with a unique associated
critical circuit $(3\to 3)$, and $C_2=\{3\}$.
(Recall our convention that Schur
complements inherit their indices from the matrices from
which they are defined, so that
$(A_2)_{33}=2$ is the top left entry of $A_2$.)
We have
$A_3=\sch(C_2,\alpha_2, A_2)
= 5\oplus 0^*4=4$,
hence, $\alpha_3=4$, with a unique associated critical circuit,
$4\to 4$, and $C_3=\{4\}$. 

To determine the critical graphs $G^c_i(A)$,
we use Proposition~\ref{lem-3}, which shows that
$G^c_\ell(A)=G^c(\hat A_\ell)$.
We already computed $G^c_1(A)=G^c(A)$.
Since $D_2=\diag(0,0,2,2)$,
and
\[
\hat A_2 =D_2^{-1}A = 
\left[\begin{array}{ccccc}
\infty & 0 & \infty & \infty\\
0 & \infty & 1  & \infty\\
-1 & \infty &\infty & 0 \\
\infty & \infty &2  & 3 
\end{array}\right]
\]
we deduce that $G^c_2(A)\setminus  G^c_1(A)$
consists of the circuit $(1\to 2\to 3\to 1)$.
Finally, $D_4=\diag(0,0,2,4)$ and
\begin{align}
\hat A_3 =D_3^{-1}A = 
\left[\begin{array}{ccccc}
\infty & 0 & \infty & \infty\\
0 & \infty & 1  & \infty\\
-1 & \infty &\infty & 0 \\
\infty & \infty &0  & 1 
\end{array}\right]
\label{e-ex-critic3}
\end{align}
which shows that $G^c_3(A)\setminus  G^c_2(A)$
consists of the circuit $(3\to 4\to 3)$.
The critical graphs are represented as follows
\begin{center}
\begin{tabular}[c]{c}
\input{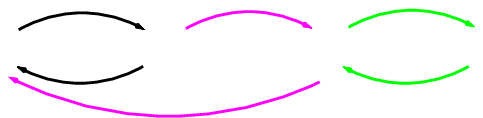}
\end{tabular}
\end{center}
Here, the graphs $G^c_{\ell}(A)$, for ${\ell}=1,2,3$ are represented in
black, magenta (medium gray), and green (light grey),
%{COLOR/BW TO FIX IN THE FINAL VERSION}
respectively; for readability, a node or arc is drawn with the
colour of the minimal graph $G^c_{\ell}(A)$ to which
it belongs.

For such a small example, the critical circuits could
be obtained by mere inspection. In general,
$G^c_\ell(A)=G^c(\hat A_\ell)$
can be computed in polynomial time thanks to Proposition~\ref{survect},
which shows that $G^c(\hat A_\ell)$ coincides with
the union of the strongly connected components
of $\sat(\hat A_\ell,V)$, for any eigenvector $V$
of $\hat A_\ell$. 
\end{example}
\subsection{Majorisation inequalities for critical values}
We now state a second majorisation result,
which should be compared with Theorem~\ref{th-gen=}.
\begin{theorem}\label{th-gam-bet}
Consider an irreducible matrix $A\in \rmin^{n\times n}$.
Let  $\Gamma=(\gamma_1\leq\cdots \leq \gamma_n)$ be the sequence of roots 
of the min-plus characteristic polynomial of $A$ and let $\beta=(\beta_1\leq
\cdots\leq\beta_n)$ be the sequence of critical values of $A$, repeated
with multiplicities. Then,
\begin{align}\label{gwb}
\Gamma \weakm \beta \enspace .
\end{align}
\end{theorem}
\begin{proof}
Let $P=\perm(\iY I\oplus A)$ be the min-plus characteristic polynomial
of $A$, $\Gamma=\corn{P}$ and $Q=(\iY\oplus \beta_1)\cdots (\iY\oplus\beta_n)$.
Let $V$ be an eigenvector of $\ainf$ (for instance any column
$\ainf^*_{\cdot j}$, since by Proposition~\ref{lem-3}, $\rho(\ainf)=\unit$
and all the nodes of $\{1,\ldots ,n\}$ are critical).
Let $W=\diag{V}$.
Since $\ainf V=V$, we get $W^{-1} \ainf W\unit=\unit$,
where $\unit$ is the vector with all entries equal to $\unit$.
Therefore, $W^{-1} A W\unit=\dinf \unit$, thus
$(W^{-1} A W)_{ij}\geq \beta_i$ for all $i,j=1,\ldots , n$.
Using~\eqref{trmin}, we get
\begin{align}
\label{tracegeqbeta}
\tr_k(A)=\tr_k(W^{-1} AW) \geq \beta_1\cdots
\beta_k \enspace .
\end{align}
Then, 
\begin{align*}
P= \iY^n\oplus \tr_1(A) \iY^{n-1}\oplus\cdots \oplus \tr_n(A)
&\geq \iY^n\oplus \beta_1 \iY^{n-1}\oplus \cdots \oplus \beta_1\cdots
\beta_n \iY^0\\
&=(\iY\oplus \beta_1)\cdots (\iY\oplus \beta_n)=Q\enspace .
\end{align*}
From Lemma~\ref{minpoly-maj}, we deduce $\corn{P} \weakm  \corn{Q}$ and
since $\Gamma=\corn{P}$ and $\corn{Q}=\beta$,
we obtain~\eqref{gwb}.
\end{proof}

We next characterise the cases where the equality
holds in~\eqref{gwb}.
We say that a graph $G$ has a \new{disjoint circuit cover}
if there is a disjoint union of circuits containing
all the nodes of $G$.
This property, which is equivalent
to the adjacency matrix of $G$ having full
\new{term rank}~\cite[Section~1.2]{brualdi91},
can be easily checked: it reduces to find
a perfect matching (or to compute a matching
of maximal cardinality) in a bipartite graph. 
\begin{theorem}\label{termrk}
Consider an irreducible matrix $A\in \rmin^{n\times n}$.
Let  $\Gamma=(\gamma_1\leq\cdots \leq \gamma_n)$ be the sequence of roots 
of the min-plus characteristic polynomial of $A$,
and let $\beta=(\beta_1\leq \cdots \leq \beta_n)$ be the sequence of
critical values of $A$ repeated with multiplicities.
For all $\ell\in\{1,\ldots, k\}$, 
where $k$ is the number of critical values of $A$,
the following assertions are equivalent:
\begin{enumerate}
\item \label{termrk.1}
$\gamma_j=\beta_j$ for $j\in\{\card{C^{\ell-1}}+1,\ldots , \card{C^\ell}\}$,
and
$\gamma_1\cdots\gamma_{\card{C^{\ell -1}}}=\beta_1\cdots\beta_{\card{C^{\ell -1}}}$;
\item\label{termrk.2}  $G_{\ell-1}^c(A)$ and
$G_{\ell}^c(A)$ have a disjoint circuit cover.
\end{enumerate}
\end{theorem}
In Theorem~\ref{termrk}, we use the convention that 
$G_0^c(A)$ is the empty graph and that it has  a disjoint circuit cover.
Recall also that $C^0=\emptyset$.

The proof of Theorem~\ref{termrk} relies
on the following lemma.
\begin{lemma}\label{equalcase}
The equality
\begin{align}
\tr_{\card C^\ell}(A)=\beta_1\cdots \beta_{\card C^\ell}
\label{termrk.0ter}
\end{align}
holds if, and only if, $G^c_\ell(A)$ has a disjoint circuit cover.
\end{lemma}
\begin{proof}
Let us first assume
that $G^c_\ell(A)$ has a disjoint circuit cover. 
Since by Proposition~\ref{lem-3},
the set of nodes of $G_\ell^c(A)$ is $C^\ell$,
there exists disjoint
elementary circuits $c_1,\ldots, c_q$ in $G_\ell^c(A)$ 
which cover all the nodes of $C^\ell$.
Let $\sigma$ be the permutation of the nodes
of $C^\ell$ which consists
of the circuits $c_1,\ldots, c_q$.
We obtain, using~\eqref{trmin}:
\begin{align*}
\tr_{\card{C^\ell}} (A)&\leq \bigotimes_{j\in C^\ell} A_{j\sigma(j)}
=\beta_1\cdots\beta_{\card{C^\ell}} \bigotimes_{j\in C^\ell} (\hat{A}_\ell)_{j\sigma(j)} = \beta_1\cdots\beta_{\card{C^\ell}} 
\enspace ,\end{align*}
since, by Proposition~\ref{lem-3},
$c_1,\ldots, c_q$ are critical circuits of $\hat{A}_\ell$ and 
$\rho(\hat{A}_\ell)=\unit$. 
Since it follows from~\eqref{tracegeqbeta}
that $\tr_{\card{C^\ell}}(A)\geq \beta_1\cdots\beta_{\card{C^\ell}}$,
we have proved~\eqref{termrk.0ter}.

Conversely, let us assume that~\eqref{termrk.0ter} holds.
Let $W$ be as in the proof of Theorem~\ref{th-gam-bet}.
By~\eqref{trmin}, there exists disjoint circuits $c_1,\ldots, c_q$
of $G(A)$ such that $|c_1|+\cdots +|c_q|=\card{C^\ell}$
and $\tr_{\card{C^\ell}}(A)=
\bigotimes_{j\in c_1\cup\cdots\cup c_q} A_{j \sigma(j)}
=\bigotimes_{j\in c_1\cup\cdots\cup c_q} (W^{-1} A W)_{j \sigma(j)}$
where $\sigma$ is the permutation of the nodes of $C^\ell$
consisting of the circuits $c_1$,\ldots, $c_q$.
Since $W^{-1} A W\unit \geq \dinf \unit$, we obtain that
$\tr_{\card{C^\ell}}(A) \geq
\bigotimes_{j\in c_1\cup\cdots\cup c_q} D_{jj}$.
If $c_1\cup\cdots\cup c_q\neq C^\ell$,
we obtain, using  $\beta_n\geq\cdots \geq
\beta_{\card{C^\ell}+1} > \beta_{\card{C^\ell}}\geq\cdots\geq 
\beta_1$, that
$\tr_{\card{C^\ell}}(A) > \beta_1\cdots \beta_{\card{C^\ell}}$,
a contradiction.
Therefore,  $c_1\cup\cdots\cup c_q=C^\ell$, and since 
$\beta_1\cdots\beta_{\card{C^\ell}}=
\tr_{\card{C^\ell}}(A)=
\bigotimes_{j\in C^\ell } A_{j \sigma(j)}$,
we get $\bigotimes_{j\in c_1\cup\cdots\cup c_q}
 (\hat{A}_\ell)_{j \sigma(j)}=\unit$.
Since $\rho(\hat{A}_\ell)=\unit$, 
the circuits $c_1,\ldots, c_q$, which are critical for $\hat{A}_\ell$, 
are critical circuits of $G_\ell^c(A)$ (by Proposition~\ref{lem-3}).
Hence, $G_\ell^c(A)$ has a disjoint circuit cover.
\end{proof}
\begin{proof}[Proof of Theorem~\ref{termrk}]
Let $P=\perm(\iY I\oplus A)$ be the min-plus characteristic polynomial
of $A$. By Lemma~\ref{carac-ci}, we have
\begin{align}
\divex{P}_{n-i}=\gamma_1\cdots \gamma_i
\leq P_{n-i}=\tr_i(A)
\enspace,
\mrm{ with equality when }
\gamma_i<\gamma_{i+1} \enspace .
\label{termrk.fond}
\end{align}
We prove \ref{termrk.2}$\implies$\ref{termrk.1}.
Assume that $G_{\ell-1}^c(A)$ and $G_{\ell}^c(A)$ have
a disjoint circuit cover.
Combining the inequality in~\eqref{termrk.fond} with~\eqref{termrk.0ter},
we get $ \gamma_1\cdots\gamma_{\card{C^\ell}} \leq  
\beta_1\cdots\beta_{\card{C^\ell}}$. Similarly,
$\gamma_1\cdots\gamma_{\card{C^{\ell-1}}} \leq  
\beta_1\cdots\beta_{\card{C^{\ell-1}}}$.
Using~\eqref{gwb}, we get the reverse inequalities
\begin{align}
\gamma_1\cdots\gamma_j\geq \beta_1\cdots\beta_j
,\mrm{ for } j=1,\ldots,n
\label{termrk.0}
\end{align}
so that
\begin{align}
\gamma_1\cdots\gamma_{\card{C^\ell}} &= \beta_1\cdots\beta_{\card{C^\ell}}
\enspace,
\label{termrk.3} \\
\gamma_1\cdots\gamma_{\card{C^{\ell-1}}} &= \beta_1\cdots\beta_{\card{C^{\ell-1}}} 
\enspace .\label{termrk.4}
\end{align}
Dividing~\eqref{termrk.3} by~\eqref{termrk.4}, we get
\begin{align}
\gamma_{\card{C^{\ell-1}}+1}\cdots \gamma_{\card{C^\ell}}
= \beta_{\card{C^{\ell-1}}+1}\cdots \beta_{\card{C^\ell}}
= \alpha_\ell^{\card C_\ell}
\label{termrk.4ter}
\end{align}
(recall that $\card C_\ell= \card C^\ell -\card C^{\ell-1}$).
Taking $j=\card C^{\ell-1}+1$ in~\eqref{termrk.0},
and using~\eqref{termrk.4}, we get 
\begin{align*}
\gamma_{\card C^{\ell -1}+1}\geq
\beta_{\card C^{\ell -1}+1}=\alpha_\ell\enspace.
\end{align*}
Since $(\gamma_i)$ is nondecreasing,
$\gamma_j\geq \gamma_{\card C^{\ell-1}+1}\geq \alpha_\ell$
holds for all $j\in \{\card C^{\ell -1}+1,\ldots,\card C^\ell\}$,
hence, if $\gamma_j>\alpha_\ell$ for some
$j\in \{\card C^{\ell -1}+1,\ldots,\card C^\ell\}$,
we would have $\gamma_{\card C^{\ell -1}+1}
\cdots \gamma_{\card C^\ell}>\alpha_\ell^{\card C_\ell}$,
contradicting~\eqref{termrk.4ter}.
Therefore,
$\gamma_{\card{C^{\ell-1}}+1}=\cdots= \gamma_{\card{C^\ell}}=\alpha_\ell=
\beta_{\card{C^{\ell-1}}+1}=\cdots= \beta_{\card{C^\ell}}$.

We next prove \ref{termrk.1}$\implies$\ref{termrk.2}.
By assumption,~\eqref{termrk.3} and~\eqref{termrk.4} hold.
Taking $j=\card C^\ell +1$ in~\eqref{termrk.0}
and using~\eqref{termrk.3}, 
we have $\gamma_{\card C^\ell+1}\geq \beta_{\card C^\ell +1}$.
Since
$\beta_{\card C^\ell +1} >\beta_{\card C^\ell}
=\gamma_{\card C^\ell}$, 
we have $\gamma_{\card C^\ell+1}>\gamma_{\card C^\ell}$,
so the equality case in~\eqref{termrk.fond} yields
\begin{align}
\gamma_1\cdots \gamma_{\card C^\ell}=\tr_{\card C^\ell}(A) \enspace.
\label{termrk.6bis}
\end{align}
Taking now $j=\card{C^{\ell-1}-1}$ in~\eqref{termrk.0},
and using~\eqref{termrk.4}, we get 
$\beta_{\card C^{\ell -1}}\geq \gamma_{\card C^{\ell -1}}$,
hence, $\gamma_{\card C^{\ell-1}+1}=\beta_{\card C^{\ell -1}+1}
>\beta_{\card C^{\ell -1}}\geq \gamma_{\card C^{\ell -1}}$,
and the equality case in~\eqref{termrk.fond} yields
\begin{align}
\gamma_1\cdots \gamma_{\card C^{\ell-1}}=\tr_{\card C^{\ell-1}}(A) \enspace.
\label{termrk.6ter}
\end{align}
It follows from Lemma~\ref{equalcase}, 
and from~\eqref{termrk.3},~\eqref{termrk.4},%
~\eqref{termrk.6bis} and~\eqref{termrk.6ter}, that
$G^c_\ell(A)$ and $G^c_{\ell-1}(A)$ have disjoint circuits covers.
\end{proof}
\begin{corollary}\label{cor-sympa}
If $G^c_{\ell -1}(A)$ and $G^c_{\ell}(A)$ have
a disjoint circuit cover, then
$\alpha_\ell$ is a root of multiplicity 
$\card C_\ell$ of the min-plus characteristic polynomial of $A$.
\end{corollary}
\begin{proof}
Since $\gamma_j=\beta_j=\alpha_\ell$
for $j\in\{\card C^{\ell-1}+1,\ldots,\card C^{\ell}\}$,
$\alpha_\ell$ is a root of multiplicity
at least $\card C^{\ell}-\card C^{\ell -1}=\card C_\ell$
of the characteristic polynomial of $A$.
Moreover, we showed 
in the proof of ``\ref{termrk.1}$\implies$\ref{termrk.2}''
of Theorem~\ref{termrk} that $\gamma_{\card C^{\ell +1}}>\gamma_{\card C^\ell}$
and $\gamma_{\card C^{\ell-1}+1}>\gamma_{\card C^{\ell -1}}$.
Thus, $\alpha_\ell$ is a root of multiplicity
exactly $\card C_\ell$ of the characteristic
polynomial of $A$.
\end{proof}
\section{Asymptotics of eigenvalues}\label{sec-cons}
\subsection{Statement and illustration of the result}\label{sec-main}
We next show that under some non-degeneracy conditions,
the first order asymptotics of the eigenvalues of $\sA_\epsilon$
are given by the critical values of $A$.
If $G$ is any graph with set of nodes $1,\ldots,n$,
and if $b\in \C^{n\times n}$,
the matrix $b^{G}$ is defined by
\begin{align*}
(b^{G})_{ij}= \begin{cases}
b_{ij}& \mrm{if $(i,j)\in G$,}\\
0 & \mrm{otherwise.}
\end{cases}
\end{align*}
Let $G$ be either the critical graph of $\ainf$ or the saturation
graph $\sat(\ainf,V)$, for any eigenvector $V$ of $\ainf$
(since by Proposition~\ref{lem-3}, all the nodes $1,\ldots, n$ belong to the critical graph of $\ainf$, we can take for $V$ any column of $\ainf^*$).

We construct the following conventional Schur complements:
\begin{equation}
s^1=a^{G},\qquad 
s^{\ell}= \sch(C^{{\ell}-1}, s^1),\; {\ell}=2,\ldots,k
\enspace. 
\label{si}
\end{equation}
The Schur complement $s^{\ell}$ is well defined as
soon as the matrix
\begin{equation} \label{ri}
r^{\ell}= a^G_{C^{{\ell}-1},
C^{{\ell}-1}}
\end{equation}
is invertible (we
adopt the convention that $r^1$ is the empty
matrix, and is invertible).
We shall also need the following matrix:
\begin{equation} \label{ti}
t^{\ell}= s^{\ell}_{C_{\ell}C_{\ell}} \enspace .
\end{equation}
When both $s^{\ell}$ and $s^{\ell-1}$ are well defined, $t^{\ell-1}$ is
invertible and we 
can compute $s^{\ell}$ from $s^{\ell-1}$ 
thanks to~\eqref{e-ident}:
\[
s^{\ell}=\sch(C_{{\ell}-1}, s^{\ell-1}) \enspace.
\]
We say that a function of $\eps$, $f(\eps)$,
is \new{of order} $\omega(\eps^\alpha)$
if $\lim_{\eps \to 0} |f(\eps) \eps^{-\alpha}|=+\infty$.
\begin{theorem}[Generalised Lidski\u\i-Vi\v{s}ik-Ljusternik theorem]\label{th-1}
Let $s^{\ell}$, $r^{\ell}$, $t^{\ell}$, ${\ell}=1,\ldots ,k$ be constructed as
in {\rm (\ref{si},\ref{ri},\ref{ti})} with $G={G^c(\ainf)}$ or 
equivalently with
$G={\sat(\ainf,V)}$ for some eigenvector $V$ of $\ainf$.
Assume that the matrix $r^{\ell}$ is invertible
for some $1\leq {\ell}\leq k$,
and let $\lambda^{\ell}_1,\ldots,\lambda^{\ell}_{m_{\ell}}$ denote
the non-zero eigenvalues of $t^{\ell}$ (here and in the sequel,
eigenvalues are repeated with multiplicities).
Then, the eigenvalues of $\sA_\eps$ can be grouped
in
\begin{enumerate}
\item $m_{\ell}$ eigenvalues with asymptotic expansions
\begin{equation}
\sL_\eps^{\ell ,j} \sim \lambda^{\ell}_j \eps^{\alpha_{\ell}} , \quad 
1\leq j\leq m_{\ell} \enspace,
\label{e-order}
\end{equation}
\item $\card{C^{{\ell}-1}}$ eigenvalues of order 
$\omega(\eps^{\alpha_{\ell}})$,
\item $\card{N^{\ell}}-m_{\ell}$ eigenvalues of order
$o(\eps^{\alpha_{\ell}})$. 
\end{enumerate}
In particular, when $t^1,\ldots, t^k$ all are invertible,
for all $1\leq \ell\leq k$, 
$\sA_\eps$ has exactly $\card{C_{\ell}}$ eigenvalues
of order $\eps^{\alpha_{\ell}}$, whose asymptotics
are given by~\eqref{e-order}.
\end{theorem}
We prove Theorem~\ref{th-1} in Section~\ref{sec-proo-th-1}.

By Proposition~\ref{survect}, the saturation graph $\sat(\ainf,V)$
(defined in Section~\ref{spectral-sec}) and the critical graph
$G^c(\ainf)$ have the same strongly connected components.
This explains why, in Theorem~\ref{th-1}, one can use
either the graph $G={G^c(\ainf)}$ or the graph $G={\sat(\ainf,V)}$.

The following result, that we also prove in  Section~\ref{sec-proo-th-1},
shows that the assumptions of the theorem are generically
satisfied, if we assume that the critical graphs
have disjoint circuit covers:
\begin{proposition}\label{prop-gen}
Let $\ell=1,\ldots, k$.
Assume that $G^c_{\ell-1}(A)$ and
$G^c_{\ell}(A)$ have disjoint circuit covers.
Then, $r^\ell$ and $t^\ell$ are generically
invertible, so that the number
of eigenvalues of $\sA_\epsilon$ with
an equivalent of the form $\lambda \epsilon^{\alpha_\ell}$,
where $\lambda\in \C\setminus\{0\}$,
is generically $\card C_\ell$.
\end{proposition}
\begin{example}\label{ex-1}
To illustrate Theorem~\ref{th-1},
consider the matrix
\begin{align}
\label{eq-mat}
\sA_\eps= \left[
\begin{array}{ccc} 
\eps & 1& \eps^4\\
0 &\eps &\eps^{-2}\\
\eps & \eps^2 & 0 
\end{array}
\right]\enspace,
\end{align}
so that $(\sA_{\eps})_{ij}\simeq a_{ij}\eps^{A_{ij}}$,
with
\[
a=
\left[\begin{array}{ccc}
1 &  1 & 1\\
1 &  1 & 1\\
1 &  1 & 1
\end{array}
\right]\enspace ,\quad
A=
\left[\begin{array}{ccc}
1 &  0 & 4\\
\infty & 1 & -2\\
1& 2 & \infty
\end{array}
\right]\enspace .
\]
We have $\rhomin(A)=-1/3$, and $\GC(A)$
consists of the critical circuit:
\begin{center}
\begin{tabular}[c]{c}
\input{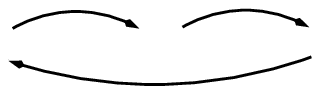}
\end{tabular}
\end{center}
so that the construction of the critical
classes stops with $C_1=\{1,2,3\}$ and $k=1$.
Then, $G^c(\ainf)=G^c_1(A)=G^c(A)$ covers all the nodes (see 
Proposition~\ref{lem-3}), hence,
\[
s^1= \left[\begin{array}{ccc}
0 &  1 & 0\\
0 & 0 & 1\\
1& 0 & 0
\end{array}
\right]\enspace .
\]
Since the spectrum of $s^1$ is $\{1,j,j^2\}$,
Theorem~\ref{th-1} shows that the spectrum
of $\sA_\eps$ consists of the three
eigenvalues
\begin{align*}
\sL^1_\eps\sim \eps^{-1/3},
\sL^2_\eps\sim j\eps^{-1/3},
\sL^3_\eps\sim j^2\eps^{-1/3}.
\end{align*}
\end{example}
\begin{example}
To give an example in which different exponents
appear, consider
\[
\sA_\eps =
\left[\begin{array}{ccccc}
\cdot & a_{12} & \cdot & \cdot\\
a_{21} & \cdot & \eps a_{23} & \cdot\\
\eps a_{31} & \cdot &\cdot & \eps^2 a_{34}\\
\cdot & \cdot &\eps^4 a_{43} &  \eps^5 a_{44}
\end{array}
\right]
\enspace ,
\]
where $a_{ij}\in \C$,
and ``$\cdot$'' denotes a zero entry.
The associated matrix of exponents $A$
is given by~\eqref{e-ex-critic}, and we saw in Example~\ref{ex-critic}
that the critical values of $A$ are $\alpha_1=0$,
$\alpha_2=2$, $\alpha_3=4$, with $C_1=\{1,2\}$,
$C_2=\{3\}$, $C_3=\{4\}$.
The critical graph $G=G^c(\hat{A})$ of the matrix $\hat{A}=\hat{A}_3$
of~\eqref{e-ex-critic3} was represented in Example~\ref{ex-critic}.
Thus, 
\[
s^1=a^{G}
=\left[\begin{array}{ccccc}
\cdot & a_{12} & \cdot & \cdot \\
a_{21} & \cdot & a_{23} & \cdot\\
a_{31}  & \cdot &\cdot &  a_{34}\\
\cdot & \cdot & a_{43} & \cdot
\end{array}
\right]
\enspace .
\]
The eigenvalues of the matrix
$t^1=\left[\begin{smallmatrix}0 & a_{12} \cr a_{21}
 & 0\end{smallmatrix}\right]$
are the square roots of $a_{12}a_{21}$.
Let us assume that $a_{12}a_{21}\neq 0$.
Then, Theorem~\ref{th-1} shows that $\sA_\eps$ has two
eigenvalues with asymptotics of the form 
$\sL_\eps\sim \xi$, where $\xi^2=a_{12}a_{21}$.
Moreover,
\[
s^2= \sch(\{1,2\}, s^1)=   
\left[\begin{array}{cc}
-a_{31}a_{21}^{-1}a_{23} &  a_{34}\\
a_{43} & \cdot
\end{array}
\right]
\enspace ,\quad t^2= -a_{31}a_{21}^{-1}a_{23} \enspace. 
\]
If we assume additionally that $a_{31}a_{23}\neq 0$,
Theorem~\ref{th-1} shows that $\sA_\eps$ has an eigenvalue
with asymptotics $\sL_\eps \sim -a_{31}a_{21}^{-1}a_{23}\eps^2$.
Finally, as soon as the matrix $r^3$ is invertible,
i.e., as soon as $\det r^3=a_{12}a_{23}a_{31}\neq 0$, 
the Schur complement
\[
t^3=s^3=
\sch(\{1,2,3\},s^1)=
a_{43}a_{23}^{-1}a_{21}a_{31}^{-1}a_{34}\enspace.
\]
is well defined. 
(When $s^2$ is well defined, that is $a_{12} a_{21}\neq 0$,
 and $a_{31}a_{23}\neq 0$,
we may obtain equivalently $t^3$ as
$\sch(\{3\}, s^2)$.)
Thus, when $a_{12}a_{23}a_{31}\neq 0$ and $a_{43}a_{21}a_{34}\neq 0$,
Theorem~\ref{th-1} shows that $\sA_\eps$ has an eigenvalue
with asymptotics $\sL_\eps \sim
a_{43}a_{23}^{-1}a_{21}a_{31}^{-1}a_{34} \eps^4$.
\end{example}
\subsection{Proof of Theorem~\ref{th-1} and Proposition~\ref{prop-gen}}
\label{sec-proo-th-1}
For the proof of Theorem~\ref{th-1}, 
we need to use the following lemma, which follows readily from
the definition of determinants. 
\begin{lemma}\label{l-circuits}
If $b,\tilde{b}\in \C^{n\times n}$ have two digraphs
$G(b)$ and $G(\tilde{b})$ whose circuits (or equivalently,
whose strongly connected components) are the same,
and if $b_{ij}=\tilde{b}_{ij}$ for all arcs $(i, j)$
belonging to circuits of $G(b)$ or $G(\tilde{b})$, then, $\det b=
\det\tilde{b}$.%\qed
\end{lemma}

Let $V$ be an eigenvector of $\ainf$ and let $\sat=\sat(\ainf,V)$.
The change of variables $\lambda=\mu \eps^{\alpha_\ell}$, for some
$1\leq \ell\leq k$, transforms
the characteristic polynomial of $\sA_\eps$ into
\begin{align*}
\det(\mu \eps^{\alpha_\ell}I -\sA_\eps )&=
\det (\eps^{D_\ell})
\det (\mu \eps^{\alpha_\ell} \eps^{D_\ell^{-1}} I 
-\eps^{D_\ell^{-1}}\sA_\eps )=
\det (\eps^{D_\ell})\sP(\eps,\mu)\nonumber\\
&\mrm{where}\; \sP(\eps , \mu) = 
\det (\mu \eps^{\alpha_\ell}\eps^{D_\ell^{-1}} I
-\eps^{D_\ell^{-1}}\eps^{\diag(V)^{-1}}
\sA_\eps \eps^{\diag(V)})
\enspace .
\end{align*}

If $C\subset L$ are finite sets,  we denote by $E_C^L$ the
$L\times L$ diagonal matrix such that 
\begin{align*}
(E_C^L)_{ii}= \begin{cases} 1 & \text{ for } i\in C\enspace,\\
0 & \text{ for } i\in L\setminus C\enspace . \end{cases}
\end{align*}
If $L=\{1,\ldots , n\}$, we shall simply write $E_C$ instead of $E_C^L$.
We have
\begin{align}
\eps^{D^{-1}}\eps^{\diag(V)^{-1}}
\sA_\eps \eps^{\diag(V)}&\tend_{\eps\to 0} a^{\sat},
& \eps^{D_\ell^{-1}}\eps^D \tend_{\eps\to 0} E_{C^\ell} \enspace,
\quad\mrm{and}\nonumber\\
\eps^{\alpha_\ell}\eps^{D_\ell^{-1}}\tend_{\eps\to 0}
E_{N^\ell}\enspace,&\label{e-diagalpha}
\end{align}
hence $\sP(\eps,\mu)\tend_{\eps\to 0} \sP(0,\mu)$, where
\begin{align*}
\sP(0,\mu)= \det (\mu E_{N^\ell}-E_{C^\ell} a^{\sat}
)\enspace .\end{align*}
Since $\sat$ and $G^c(\ainf)$ have the same strongly connected components
(by Proposition~\ref{survect}), Lemma~\ref{l-circuits} yields:
\begin{align*}
\sP(0,\mu)= \det(\mu E_{N^\ell}-E_{C^\ell} a^{G^c(\ainf)}
) \enspace .
\end{align*}
The same arguments also show that the
invertibility of the matrix $r^\ell$ is independent of the choice of
$G=\sat$  or $G=G^c(\ainf)$ in~\eqref{si}.
Hence, if $s^\ell$, $r^\ell$, $t^\ell$ are constructed 
as in {\rm (\ref{si},\ref{ri},\ref{ti})} with either $G={\sat}$ or $G=G^c(\ainf)$,
and if $r^{\ell}$ is invertible, then
\begin{align*}
\sP(0,\mu)&= \mu^{\card{N^{\ell +1}}}
\det (\mu  E_{C_\ell}^{C^\ell}-a^{G}_{C^\ell,C^\ell})
= \mu^{\card{N^{\ell +1}}} \det (-r^\ell) \det(\mu I -t^\ell ) \enspace . 
\end{align*}
From Lemma~\ref{cont-roots}
applied to $\sP(\eps,\iY)$, there exists $\card{N^\ell}$
continuous functions $\eps\mapsto \sL_\eps^{m,j}$, 
with $j=1,\ldots, \card{C_m}$ and $m=\ell,\ldots, k$,
such that $\sL_\eps^{m,j}$ are the roots of $\sP(\eps,\iY)$
for all $\eps$ small enough. Hence, 
$\sL_0^{\ell,j}$ are the eigenvalues of $t^\ell$ and 
$\sL_0^{m,j}=0$ for $m>\ell$. The other roots of $\sP(\eps,\iY)$
tend to infinity.
This shows Theorem~\ref{th-1}.
%\qed

We finally prove Proposition~\ref{prop-gen}.
If the set of nodes of $G^c_{\ell-1}(A)$ can be covered
by disjoint circuits, it follows from Proposition~\ref{lem-3}
that these circuits also belong to 
$G^c(\hat A)\cap C^{\ell-1}\times C^{\ell-1}$. By definition
of $r^\ell$, for generic values of $a=(a_{ij})$,
these circuits belong 
to the graph of $r^\ell$, which implies
that the determinant of $r^\ell$
is generically non-zero. Thus, $r^\ell$ is generically invertible.
The same argument shows that if $G^c_{\ell}(A)$ can
be covered by disjoint circuits,
$r^{\ell+1}$ is generically
invertible, and since $t^\ell=\sch(C^{\ell-1},r^{\ell+1})$
is the Schur complement of the generically invertible
$C^{\ell -1}\times C^{\ell -1}$ submatrix of $r^{\ell+1}$,
namely $r^\ell$, in the generically invertible matrix $r^{\ell+1}$, 
$t^\ell$ must also be generically invertible.
Thus, $m_\ell=\card C_\ell$ generically in Theorem~\ref{th-1}.
%\qed

\section{Asymptotics of eigenvectors}\label{sec-vec}
\subsection{Statement and illustration of the result}
We now consider eigenvectors.
\begin{theorem}\label{th-2}
Let $s^\ell$, $r^\ell$, $t^\ell$, $\ell=1,\ldots ,k$ be constructed as
in Theorem~\ref{th-1}.
Assume that the matrix $r^{\ell}$ is invertible, for some $1\leq \ell \leq k$,
that $\mu\neq 0$ is a simple eigenvalue of $t^\ell$,
and let $V$ be any eigenvector of $\hat{A}_\ell$.
Then, the equation 
\begin{align}
(\mu E_{N^\ell}- a^{\sat(\hat A_\ell,V)}) w=0 ,
\label{e-d-w}
\end{align}
has a unique solution $w=(w_j)\in \C^n\setminus\{0\}$
up to a multiplicative constant.
Moreover, there is a unique eigenvalue
$\sL_\eps$ with asymptotics
$\sL_\eps \sim \mu \eps^{\alpha_\ell}$,
and if $w_i\neq 0$, any eigenvector 
$\sV_\eps$ associated to this $\sL_\eps$ satisfies $(\sV_\eps)_i\not =0$
for $\eps$ small enough, and
\begin{align}
\label{e-ratio}
\frac{(\sV_\eps)_j}{(\sV_\eps)_i}
\simeq \frac{w_j \eps^{ V_j}}{w_i \eps^{V_i}}
\enspace ,
\mrm{ for } j\in \{1,\ldots,n\} \enspace .
\end{align}
\end{theorem}
We prove Theorem~\ref{th-2} in Section~\ref{proof-th2}.

\begin{example}\label{ex-1-eigenvect}
To illustrate Theorem~\ref{th-2},
let us pursue the analysis of Example~\ref{ex-1}.
We already showed that the eigenvalues
of the matrix~\eqref{eq-mat} have asymptotic
equivalents of the form $\xi \epsilon^{-1/3}$,
where $\xi$ is a cubic root of $1$.
When $\mu=\xi$, any solution of~\eqref{e-d-w} (with $\ell=1$),
is proportional to $w=[1,\xi,\xi^2]^T$.
Since $A$ has a unique critical class,
$C_1=\{1,2,3\}$, 
by Theorem~\ref{min-plus-spectral2},
$A$ has a unique
eigenvector, up to a scalar factor,
and we can take $V =[0,- 1/3, 4/3]^T=\ainf^*_{\cdot, 1}$.
Theorem~\ref{th-2} shows that
any eigenvector $\sV_\eps$ associated to
the eigenvalue $\xi \eps^{-1/3}$ is 
equivalent to
\[
[1, \xi \eps^{-1/3}, \xi^2\eps^{4/3}]^T \enspace ,
\]
up to a scalar factor.
\end{example}
When $w_j=0$, Theorem~\ref{th-2}
gives a poor information on the asymptotics of $(\sV_\eps)_j$.
Moreover, when $\hat{A}_\ell$ has several critical
classes (so that the eigenvector $V$ is non unique)
the non-zero character of $w_j$ depends
in a critical way of the eigenvector $V$ which is
selected. 
\begin{example}\label{ex-canonical}
The following example illustrates the importance
of the choice of the eigenvector $V$ in Theorem~\ref{th-2}.
Consider
\[
\sA_\eps= \left[\begin{array}{ccc}
1 & \eps & \eps^3\\
-2\eps & \eps^2& \cdot \\
\eps^3 & \cdot & 2\eps^2
\end{array}
\right]
\]
which is such that 
$(\sA_\eps)_{ij}\simeq a_{ij}\eps^{A_{ij}}$
with
\[
A=\left[\begin{array}{ccc}
0 & 1 & 3\\
1 & 2&\infty \\
3 & \infty & 2
\end{array}
\right],
\;
a= \left[\begin{array}{ccc}
1 & 1 & 1\\
-2 & 1&\cdot \\
1 & \cdot & 2
\end{array}\right]\enspace .
\]
We have $\alpha_1=\rhomin(A)=0$, with a unique critical circuit $(1\to 1)$.
Hence, $C_1=\{1\}$, and 
\[
A_2=\sch(\{1\},A)=  
\left[\begin{array}{cc}
2&4 \\
4 & 2
\end{array}
\right]\enspace .
\]
Thus, $\alpha_2=2$, $C_2=\{2,3\}$. We have
\[
\ainf= 
\left[\begin{array}{ccc}
0 & 1 & 3\\
-1 & 0&\infty \\
1 & \infty & 0
\end{array}
\right].
\]
Since the critical graph of $\ainf$,
which is the union of the complete graph on $\{1,2\}$,
and of the loop $(3\to 3)$, has
two strongly connected components,
$\{1,2\}$, and $\{3\}$, the eigenspace of $\ainf$
is spanned by the two vectors $\ainf^*_{\cdot, i}$,
$i=1,3$. Let us take 
\[V=\ainf^*_{\cdot, 3}=[3,2,0]^T
\enspace,
\]
for which the saturation graph is obtained by adding  the arc
$(1\to 3)$ to the critical graph of $\ainf$. 
Taking $G=\sat(\ainf ,V)$ in~\eqref{si}, we get
\begin{align}
s^1=a^{\sat(\ainf ,V)}= \left[\begin{array}{ccc}
1 & 1 & 1\\
-2 & 1&\cdot \\
\cdot & \cdot & 2
\end{array}\right]\enspace .
\label{e-e3}
\end{align}
Since $t^1=1$, Theorem~\ref{th-1} shows that
$\sA_\eps$ has a root with asymptotics
$\sL_\eps \sim 1$, and since $t^2=s^2=\sch(\{1\},s^1)
= \left[\begin{smallmatrix} 3& 2\cr 0 & 2\end{smallmatrix}\right]$
has roots $2,3$, $\sA_\eps$ has two eigenvalues with respective asymptotics
$\sL_\eps\sim 2\eps^2$, and $\sL_\eps\sim 3\eps^2$.
Let us compute for instance the asymptotics
of the eigenvector $\sV_\eps$ associated
to $\sL_\eps\sim 2\eps^2$, using Theorem~\ref{th-2}
with $\ell=2$ and $\mu=2$ (thus $\hat{A}_2=\ainf$).
With the previous choice of $V$,
we need to solve the system~\eqref{e-d-w}, which,
by~\eqref{e-e3}, specialises to
\[
w_1 + w_2 + w_3 =0,\; -2w_1 - w_2 =0,\; 0 =0 \enspace.
\]
All the solutions of this system are proportional to $w=[1,-2, 1]^T$. Thus,
Theorem~\ref{th-2} shows that up to a multiplicative constant,
\[
\sV_\eps \sim [ \eps^3, -2\eps^2, 1]^T \enspace .
\]
Consider now the alternative choice of $V$:
\[
V=\ainf^*_{\cdot,1}=[0,-1,1]^T
\enspace .
\]
Then, $\sat(\hat{A}_2,V)$
is obtained by adding the arc $(3\to 1)$ to the critical graph
of $\ainf$. Theorem~\ref{th-2} yields
that $(\sV_\eps)_i \simeq w_i \eps^{V_i}$, where
\[ w_1+ w_2=0,\; -2w_1-w_2 =0 ,\; w_1 =0\]
and since all the solutions $w$ are proportional to $[0,0,1]^T$, 
we learn only from~\eqref{e-ratio} that 
$(\sV_\eps)_1/(\sV_\eps)_3 \simeq 0 \eps^{-1}$, and  
$(\sV_\eps)_2/(\sV_\eps)_3\simeq 0 \eps^{-2}$, 
a very poor information. 
\end{example}

\begin{remark}
When $\mu$ is not a simple root of $t^\ell$,
the first order asymptotics of the eigenvector
may be ruled by higher order terms in the
expansions of the entries of $\sA_\epsilon$,
see~\cite{ABG96} for a special case.
\end{remark}
\subsection{Proof of Theorem~\ref{th-2}}
\label{proof-th2}
We first observe that by Theorem~\ref{th-1},
there is only one eigenvalue
$\sL_\epsilon$ of $\sA_\epsilon$
equivalent to $\mu\epsilon^{\alpha_\ell}$.
Then the associated eigenvector, $\sV_\epsilon$,
is unique, up to a multiplicative
constant, since for $\epsilon$ small enough,
$\sL_\epsilon$ is a simple eigenvalue of $\sA_\epsilon$.

To prove Theorem~\ref{th-2},
we perform the change of variables
$\sV_\eps=\eps^{\diag V} \sW_\eps$ and
$\sL_\eps=\sM_\eps \eps^{\alpha_\ell}$, where $\sM_\eps\to \mu$
when $\eps\to 0$.
After multiplying $\sV_\eps$ by a constant, we 
may assume that $\sum_{1\leq j\leq n} |(\sW_\eps)_j|=1$. 
From $\sA_\eps\sV_\eps=\sL_\eps \sV_\eps$, we get
\[
\eps^{D_\ell^{-1}}\eps^{(\diag V)^{-1}}\sA_\eps \eps^{\diag V}
\sW_\eps = \sM_\eps \eps^{D_\ell^{-1}}\eps^{\alpha_\ell} \sW_\eps \enspace,
\]
where $\eps^{D_\ell^{-1}}\eps^{\diag(V)^{-1}}
\sA_\eps \eps^{\diag V}
\to 
a^{\sat(\hat A_\ell, V)}$
when $\eps\to 0$. Together with~\eqref{e-diagalpha},
this implies that any limit point $w$ of $\sW_\eps$ when $\eps\to 0$
satisfies 
\begin{align}
a^{\sat(\hat A_\ell, V)} w = \mu E_{N^\ell} w
,\quad \mrm{and}\;
|w_1|+\cdots+|w_n|=1
\enspace  .
\label{e-lp}
\end{align}
To show that the solution $w$ of~\eqref{e-lp} is unique,
up to the multiplication by a complex number of modulus $1$,
we shall prove that $\mu E_{N^\ell}-a^{\sat(\hat A_\ell, V)}$
has rank $n-1$. 

Since, by Proposition~\ref{survect},
 $\sat(\hat A_\ell ,V)$ and $G^c(\hat A_\ell)=G^c_\ell(A)$ 
have the same 
strongly connected components, applying Lemma~\ref{l-circuits} to 
the matrices $b=b(\lambda) =\lambda  E_{N^\ell}-
a^{\sat(\hat A_\ell,V)}$
and $\tilde{b}=\tilde{b}(\lambda)=\lambda E_{N^\ell}-a^{G^c_\ell(A)}$,
with $\lambda\in \C$,
we get $\det b(\lambda)=\det \tilde{b}(\lambda)$.
Moreover, since $G^c_\ell(A)$ and the restriction of $G^c(\ainf)$ to
$C^\ell$ have the same strongly connected components
(see Proposition~\ref{lem-3}), then
by Lemma~\ref{l-circuits} again, $\det \tilde{b}(\lambda)= 
\det (\lambda  E_{C_\ell}^{C^\ell}-r^{\ell+1})
\lambda^{\card{N^{\ell +1}}}=\det(-r^{\ell})
\det (\lambda I -t^\ell)\lambda^{\card{N^{\ell +1}}}$,
which yields:
\begin{align}\label{blambda}
 \det b(\lambda)=\det(-r^{\ell})
\det (\lambda I-t^\ell)\lambda^{\card{N^{\ell +1}}}\enspace .
\end{align}
Hence, $\det b(\mu)=0$ since $\mu$ is an eigenvalue of $t^\ell$, and
$\mu E_{N^\ell}-a^{\sat(\hat A_\ell, V)}$
has rank $<n$. Since $\mu$ is a simple eigenvalue of $t^\ell$ and 
$\mu\not = 0$,
$\mu$ is a simple root of the equation $\det b(\lambda)=0$. Hence, 
the partial derivative $\partial_\lambda\det b(\lambda)$, evaluated
at $\lambda=\mu$, is non-zero, which implies that 
there is a subset $L$ of $\{1,\ldots , n\}$,
of cardinality $n-1$, such that
$\det (b(\mu)_{L,L})\not =0$, which shows that 
$\mu E_{N^\ell}-a^{\sat(\hat A_\ell,V)}$ has 
rank $n-1$.  Thus, \eqref{e-d-w} has only
one non-zero solution, up to a scalar multiple,
which implies that all the solutions of~\eqref{e-lp}
are of the form $\zeta w$, where $\zeta\in\C$ is such that
$|\zeta|=1$, and $w$ is any solution
of~\eqref{e-lp}. Let us pick $i$ such that $w_i\neq 0$.
Since all the limit
points of $\sW_\epsilon$ are of the form $\zeta w$,
with $|\zeta|=1$, we get $(\sW_\epsilon)_j/(\sW_\epsilon)_i\to
w_j/w_i$ when $\epsilon \to 0$, and since 
$\sV_\eps=\eps^{\diag V} \sW_\eps$, we get~\eqref{e-ratio}.
\subsection{On the choice of the eigenvector $V$}
\label{rk-formal}
We now show that there is, in some sense, a
canonical choice of $V$ in Theorem~\ref{th-2}.
Denote by $C^\ell_1,\ldots, C^\ell_{\nu_\ell}$ the critical classes
of  $\hat{A}_\ell$, 
and by $C_\ell^1,\ldots, C_\ell^{\nu_\ell}$ their restrictions to
$C_\ell$. By Proposition~\ref{lem-3},
$C^\ell_1,\ldots, C^\ell_{\nu_\ell}$ are the strongly
connected components of $G^c(\ainf )\cap C^\ell\times C^\ell$ and they
cover $C^\ell$.
Moreover, one can deduce from Proposition~\ref{lem-alpha}, that
for $\nu =1,\ldots,\nu_\ell$, $C_\ell^\nu$ is either the empty set or
a critical class of the matrix ${A}_\ell$, and that
$C_\ell^1\cup\cdots\cup C_\ell^{\nu_\ell}=C_\ell$.
Then, when $r^\ell$ is invertible,
the characteristic polynomial of $t^\ell$ can be factored as 
\begin{align}
\det(\lambda I - t^\ell)=
Q_\ell^1(\lambda)\cdots Q_\ell^{\nu_\ell}(\lambda)\label{lem-factored}
\end{align}
where $Q_\ell^\nu (\lambda)=\det(\lambda I -t^\ell_{C_\ell^\nu ,
C_\ell^\nu } )$ if $C_\ell^\nu\neq \emptyset$ and
$Q_\ell^\nu (\lambda)=1$ otherwise.
Indeed, taking $G=G^c(\ainf )$ in~\eqref{si}, using the fact that
 $C^\ell_1,\ldots, C^\ell_{\nu_\ell}$ are the strongly
connected components of $G^c(\ainf )\cap C^\ell\times C^\ell$, 
and using the block triangular structure of $\lambda E_{N^\ell}-a^{G^c(\ainf )}$,
we get
\begin{align}
 \det(-r^\ell)  \det(\lambda I - t^\ell) &=
\det (\lambda E_{C_\ell}^{C^\ell} -a^{G^c(\ainf )}_{C^\ell C^\ell})
\label{factorisation}\\
&= \prod_{\nu=1}^{\nu_\ell} 
\det (\lambda E_{C_\ell^\nu}^{C^\ell_\nu} -a^{G^c(\ainf )}_{C^\ell_\nu
 C^\ell_\nu}) \nonumber\\
&= \det(-r^\ell) \prod_{\nu=1,\ldots, \nu_\ell,\; C_\ell^\nu\neq\emptyset }
\det (\lambda I -t^\ell_{C_\ell^\nu  C_\ell^\nu}) 
\enspace .\nonumber
\end{align}
Since $r^\ell$ is invertible, this shows~\eqref{lem-factored}.
Thus, if $\mu\neq 0$ is a simple root of $\det(\lambda I-t^\ell)$,
there is a unique $\nu\in\{1,\ldots ,\nu_\ell\}$ such that $\mu$ is
a root of the polynomial $Q_\ell^\nu(\lambda )$.
Denote by $\nu(\mu)$ this index.
Let $V$ be an eigenvector of $\hat{A}_\ell$, for instance
$V=(\hat A_\ell)^*_{\cdot, j}$ with $j\in C^\ell$.
By the same arguments as in the proof of~\eqref{blambda},
one can show that~\eqref{factorisation} remains valid if we
replace $G^c(\hat A_\ell)$ by $\sat(\hat A_\ell, V)$.
Hence, for any $\nu\neq \nu(\mu)$,
$(\mu E_{N^\ell}-a^{\sat(\hat A_\ell, V)})_{C^\ell_\nu,C^\ell_\nu}$
is invertible. Moreover, since $\mu\neq 0$, 
$(\mu E_{N^\ell }-a^{\sat(\hat A_\ell, V)})_{N^{\ell+1},N^{\ell+1}}$
is invertible. 
One can then deduce, using the block triangular structure of
$\mu E_{N^\ell}-a^{\sat(\hat A_\ell, V)}$,  that if 
there is no path from $i$ to $\Cellmuh$ in 
$\sat(\hat A_\ell, V)$, then
$w_i=0$. In particular, using Proposition~\ref{form-sat}, one
deduce that if $V= (\hat A_\ell)^*_{\cdot, j} $ with $j\in C^\ell\setminus
\Cellmuh$, then there exists a final class $C^\ell_\nu$ of 
$\sat(\hat A_\ell, V)$ different from $\Cellmuh$, hence
$w_i=0$ for all $i\in C^\ell_\nu$.
This observation explains Example~\ref{ex-canonical},
and it also suggests that the choice  $V= (\hat A_\ell)^*_{\cdot, j} $
with $j\in \Cellmuh$ is canonical (note that different
choices of $j\in \Cellmuh$ yield proportional vectors $V$).
However, in the case of eigenvectors, there does
not seem to be a simple analogue
of Proposition~\ref{prop-gen} (characterising
the cases where generically $w$ has non-zero entries).

\section{The theorem of Vi\v sik, Ljusternik, and Lidski\u\i\ revisited}
\subsection{Statement of the theorem}
We now show that the theorem of Vi\v sik and Ljusternik~\cite{vishik}
and Lidski\u\i~\cite{lidskii} 
can be obtained as a corollary of Theorem~\ref{th-1},
and that Theorem~\ref{th-1} allows to solve cases
to which the classical result does not apply.
The presentation of this subsection is inspired by~\cite{mbo97},
that the reader may consult for a general discussion
of the theory of Vi\v sik, Ljusternik, and Lidski\u\i.

Lidski\u\i~\cite{lidskii} considers a matrix of the form $\sA_\eps
= \sA_0+\eps b$, where $b\in \C^{n\times n}$
and $\sA_0\in\C^{n\times n}$ is a nilpotent matrix.
We shall need specific notations for Jordan matrices.
Let $N[q]$ denote the $q\times q$ nilpotent matrix
such that $(N[q])_{i,j}=1$ if $j=i+1$,
and $(N[q])_{i,j}=0$ otherwise. For $m\geq 1$, 
we define $\Nil(m:q)
= N(q)\d+\cdots \d+ N(q)$ ($m$-times), 
where $\d+$ denotes the block diagonal sum, and, given
a decreasing sequence $q_1>q_2>\ldots >q_k\geq 1$,
and $m_1,\ldots,m_k\geq 1$, we define
$\Nil(m_1,\ldots,m_k:q_1,\ldots,q_k)
=\Nil(m_1:q_1)\d+\cdots\d+\Nil(m_k:q_k)$.
For instance, when $q_1=3,m_1=1,q_2=2,m_2=2,q_3=1,m_3=1$, we
have
\def\dbox{\Box\!\!\!\cdot}
\begin{align}
\Nil(1,2,1:3,2,1)=\left[
\begin{array}{ccc|cc|cc|c}
\cdot & 1&\cdot&\cdot&\cdot  &\cdot &\cdot&\cdot\\
\cdot &\cdot &1&\cdot  &\cdot &\cdot&\cdot\\
\dbox &\cdot &\cdot &\dbox &\cdot &\dbox &\cdot&\dbox\\\hline
\cdot &\cdot &\cdot &\cdot &1 &\cdot &\cdot&\cdot\\
\dbox &\cdot &\cdot &\dbox &\cdot &\dbox &\cdot&\dbox\\\hline
\cdot &\cdot &\cdot &\cdot &\cdot &\cdot &1&\cdot\\
\dbox &\cdot &\cdot &\dbox &\cdot &\dbox &\cdot&\dbox\\\hline
\dbox &\cdot &\cdot &\dbox &\cdot &\dbox &\cdot&\dbox
\end{array}\right]\enspace,
\label{ex-new2}
\end{align}
where, again, ``$\cdot$'' represents $0$
(why some zero entries are written inside boxes will be explained below).
We consider the case where $\sA_0$ is equal
to $\Nil(m_1,\ldots,m_k:q_1,\ldots,q_k)$.
If $1\leq \ell\leq k$, we finally 
define the $(m_1+\cdots+m_\ell)\times (m_1+\cdots+m_\ell)$
submatrix $\Phi_\ell$ of $b$, obtained by 
considering only the bottom rows and first columns of the Jordan cells
of sizes $q_i\times q_j$, $i,j=1,\ldots ,\ell$.
For instance, in the case of~\eqref{ex-new2},
\[
\Phi_1 = \left[\begin{array}{c} b_{31}
\end{array}\right]\enspace, \qquad 
\Phi_2= \left[\begin{array}{ccc}
b_{31}&b_{34} & b_{36}\\
b_{51}&b_{54} & b_{56}\\
b_{71}& b_{74} & b_{76}
\end{array}\right] \enspace,
\qquad
\Phi_3= \left[\begin{array}{cccc}
b_{31}&b_{34} & b_{36}&b_{38}\\
b_{51}&b_{54} & b_{56}&b_{58}\\
b_{71}& b_{74} & b_{76}&b_{78}\\
b_{81}& b_{84} & b_{86}&b_{78}\\
\end{array}\right] \enspace.
\]
The corresponding positions in the matrix $\sA_0$ were depicted
by boxes in~\eqref{ex-new2}. By convention, $\Phi_0$ is the empty
matrix, and is invertible.
\begin{corollary}[{\cite[Th.~1]{lidskii}}]\label{cor-l}
Assume that both $\Phi_{\ell-1}$ and $\Phi_\ell$ are invertible, for some 
$1\leq \ell \leq k$, and let $\lambda_1,\ldots,\lambda_{m_\ell}$
denote the eigenvalues of $\sch(\Phi_{\ell-1},\Phi_\ell)$.
Then, $\sA_\eps$ has $m_\ell q_\ell$ eigenvalues
with asymptotics
\[ \sL_\eps \sim \xi \eps^{1/q_\ell},\; 
\mrm{where}\; \xi^{q_\ell}= \lambda_{i}\;\mrm{and}\; i=1,\ldots,m_\ell
\]
(for each $\lambda_i$, all the $q_\ell$-th roots $\xi$ of $\lambda_i$
are taken). \end{corollary}

Of course, Corollary~\ref{cor-l} can be stated
in an equivalent ``coordinate free'' way, by using left and right
eigenvectors associated to the different Jordan blocks,
see~\cite{lidskii}.
In fact, Moro, Burke, and Overton
observed that we need not require $\Phi_\ell$ to be invertible in
Corollary~\ref{cor-l}: when $\Phi_\ell$ is singular,
\cite[Th.~2.1]{mbo97} shows that to each eigenvalue $\lambda_i\in \C$
of $\sch(\Phi_{\ell-1},\Phi_\ell)$ corresponds $q_\ell$ eigenvalues
of $\sA_\eps$ with asymptotics $ \sL_\eps=
\xi \eps^{1/q_\ell}+o(\eps^{1/q_\ell})$ where $\xi^{q_\ell}= \lambda_{i}$.

\subsection{Derivation of Corollary~\ref{cor-l}}\label{sec-cor}
Let us denote by $\Anil(m_1,\ldots,m_k:q_1,\ldots,q_k)$
the matrix of exponents associated to $\sA_\epsilon=\Nil(m_1,\ldots,m_k:q_1,\ldots,q_k)+\epsilon b$: $\Anil(m_1,\ldots,m_k:q_1,\ldots,q_k)$
is obtained from $\Nil(m_1,\ldots,m_k:q_1,\ldots,q_k)$
by exchanging zeros and ones. For instance, $\Anil(1:2)=
\left[\begin{smallmatrix} 1&0\cr 1 & 1\end{smallmatrix}\right]$
corresponds to $\sA_\eps=
\left[\begin{smallmatrix} 0&1\cr 0 & 0\end{smallmatrix}\right]
+ \epsilon b$. 
The following lemma is straightforward.
\begin{lemma}\label{l-rec}
Let $q_1>\cdots>q_k\geq 1$, $m_1,\ldots,m_k\geq 1$.
The matrix $\Anil(m_1,\ldots,m_k:q_1,\ldots,q_k)$ has
min-plus eigenvalue $1/q_1$, and set of critical
nodes $C_1=\{1,\ldots,m_1q_1\}$. Moreover,
\begin{align}
\sch(C_1,1/q_1, \Anil(m_1,\ldots,m_k:q_1,\ldots,q_k))
= \Anil(m_2\ldots,m_k:q_2,\ldots,q_k)
\enspace.
\label{e-rec}
\end{align}
\end{lemma}
It follows from Lemma~\ref{l-rec} and in particular,
from the recursive property~\eqref{e-rec}, that the sequence
of critical values of $\Anil(m_1,\ldots,m_k:q_1,\ldots,q_k)$
is $(\alpha_1,\ldots, \alpha_k)=(1/q_1,\ldots, 1/q_k)$, and
that the associated critical classes 
are $C_1=\{1,\ldots,m_1q_1\}$,
\ldots, $C_k=\{\sum_{\ell=1}^{k-1} m_{\ell}q_{\ell}+1,\ldots, 
\sum_{\ell=1}^{k} m_\ell q_\ell \}$.
Recall that the diagonal matrix $D$ is defined
from the $\alpha_\ell$ and $C_\ell$.

An eigenvector $\Vnil(m_1,\ldots,m_k:q_1,\ldots,q_k)$
of $D^{-1}\Anil(m_1,\ldots,m_k:q_1,\ldots,q_k)$
can be built as follows. For all $q \geq 1$,
we set $\Vnil(1:q)= [0, 1/q,\ldots, (q-1)/q]^T$, 
then, for $m\geq 1$, we define
$\Vnil(m:q)=\Vnil(1:q)\d+\cdots\d+\Vnil(1:q)$
($m$-times), where $\d+$ denotes the concatenation
of vectors, and, finally, we set
$\Vnil(m_1,\ldots,m_k:q_1,\ldots,q_k)=
\Vnil(m_1:q_1)\d+\cdots\d+\Vnil(m_k:q_k)$.
It is easy to see that $V=\Vnil(m_1,\ldots,m_k:q_1,\ldots,q_k)$
is an eigenvector of $\ainf =D^{-1}\Anil(m_1,\ldots,m_k:q_1,\ldots,q_k)$,
and that the corresponding saturation graph is the union
of the graph of $\Nil(m_1,\ldots,m_k:q_1,\ldots,q_k)$
and of the arcs $(i,j)$, where $i$ is the index of 
a bottom row of a Jordan block of $\Nil(m_1,\ldots,m_k:q_1,\ldots,q_k)$,
and $j$ is the index of the left column 
of a Jordan block of $\Nil(m_1,\ldots,m_k:q_1,\ldots,q_k)$.
Since $\sat(\ainf , V)$ is strongly connected, it is also equal
to $G^c(\hat A)$.
For instance, for $\sA_\eps =\Nil(1,2,1:3,2,1)+\epsilon b$,
and $G=G^c(\ainf)=\sat(\ainf , V)$, we get
\begin{align}
a^{G}= 
\left[\begin{array}{ccc|cc|cc|c}
\cdot&1&\cdot&\cdot&\cdot&\cdot&\cdot&\cdot\\
\cdot&\cdot&1&\cdot&\cdot&\cdot&\cdot&\cdot\\
b_{31}&\cdot&\cdot&b_{34}&\cdot&b_{36}&\cdot&b_{38}\\\hline
\cdot&\cdot&\cdot&\cdot&1&\cdot&\cdot&\cdot\\
b_{51}&\cdot&\cdot&b_{54}&\cdot&b_{56}&\cdot&b_{58}\\\hline
\cdot&\cdot&\cdot&\cdot&\cdot&\cdot&1&\cdot\\
b_{71}&\cdot&\cdot&b_{74}&\cdot&b_{76}&\cdot&b_{78}\\\hline
b_{81}&\cdot&\cdot&b_{84}&\cdot&b_{86}&\cdot&b_{88}\\
\end{array}\right]
\enspace.
\label{e-sat-8}
\end{align}
The statement of Corollary~\ref{cor-l} becomes a special case of 
the statement of Theorem~\ref{th-1}, provided
the following identity is proved:
\begin{align*}
\det(\lambda I-t^\ell) = 
\det(\lambda^{q_\ell }-\sch(\Phi_{\ell -1},\Phi_\ell )),\quad \ell =1,\ldots, k
\enspace.
\end{align*}
This can be seen immediately by noting that 
$t^\ell$ is a matrix of cyclicity $q_\ell$, which
can be put, by applying a transformation $t^\ell\mapsto P_\ell t^\ell
P_\ell^{-1}$,
for some permutation matrix $P_\ell$, in block circular form  
\begin{align}
P_\ell t^\ell P_\ell^{-1}=\left[\begin{array}{cc} \cdot &
 I_{m_\ell  (q_\ell -1)} \\
\sch(\Phi_{\ell -1},\Phi_\ell ) &  \cdot
\end{array}\right] \enspace ,
\label{e-tedious}
\end{align}
where $I_{q}$ is the identity matrix of order $q$,
and where the ``$\cdot$'' represent blocks with $0$ values. 

Indeed, by \eqref{ti} and \eqref{si}, we get:
\begin{align}
t^\ell = \sch(C^{\ell -1}, a^{G}_{C^\ell C^\ell})
\end{align}
and for each $\ell =1,\ldots , k$, there exists a 
matrix $Q_\ell $ corresponding to a permutation of $C^\ell$
preserving $C_\ell $, such that in block form we get:
\[
Q_\ell  a^{G}_{C^\ell,C^\ell} Q_\ell ^{-1}
=
\left[\begin{array}{cc|cc}
\cdot&I_{m_1(q_1-1)+\cdots + m_{\ell -1} (q_{\ell -1}-1)}&\cdot&\cdot\\
\Phi_\ell ^{11} &\cdot&\Phi_\ell ^{12}&\cdot\\\hline
\cdot&\cdot&\cdot&I_{m_\ell  (q_\ell -1)}\\
\Phi_\ell ^{21}&\cdot&\Phi_\ell ^{22}&\cdot
\end{array}\right], 
\]
where $\Phi_\ell =\SMALLMATRIX{ \Phi_\ell ^{11}&\Phi_\ell ^{12}\\
\Phi_\ell ^{21}&\Phi_\ell ^{22}}$ and $\Phi_\ell ^{11}=\Phi_{\ell -1}$
(for each $\ell$, the indices of $\Phi_\ell ^{22}$ correspond to the nodes
of $C_\ell =\{\sum_{i=1}^{\ell-1} m_{i}q_{i}+1,\ldots, 
\sum_{i=1}^{\ell} m_{i}q_{i}\}$ of the form $\sum_{i=1}^{\ell-1} m_{i}q_{i}
+ m q_\ell $ with $m=1,\ldots , m_\ell $). 
Hence, taking for $P_\ell $ the restriction of $Q_\ell $ to $C_\ell $, and using the fact 
 that $\SMALLMATRIX{\cdot&\Psi\\ \Phi&\cdot}^{-1}
= \SMALLMATRIX{\cdot&\Phi^{-1}\\ \Psi^{-1}&\cdot}$
 for all invertible matrices $\Psi$ and $\Phi$, we get \eqref{e-tedious}.

For instance, in the special case of~\eqref{e-sat-8}, and  $\ell =2$,
we get 
\[
Q_2a^{G}_{C^2,C^2} Q_2^{-1}
=
\left[\begin{array}{ccc|cccc}
\cdot&1&\cdot&\cdot&\cdot&\cdot&\cdot\\
\cdot&\cdot&1&\cdot&\cdot&\cdot&\cdot\\
b_{31}&\cdot&\cdot&b_{34}&b_{36}&\cdot&\cdot\\\hline
\cdot&\cdot&\cdot&\cdot&\cdot&1& \cdot\\
\cdot&\cdot&\cdot&\cdot&\cdot&\cdot&1\\
b_{51}&\cdot&\cdot&b_{54}&b_{56}&\cdot&\cdot\\
b_{71}&\cdot&\cdot&b_{74}&b_{76}&\cdot&\cdot
\end{array}\right], 
\]
and, 
\begin{align*}
P_2  t^2 P_2 ^{-1} &=
\left[\begin{array}{c|c}
\begin{array}{ccc}
\cdot & &\cdot\\\cdot&&\cdot
\end{array}
&
\begin{array}{cc}
1 & \cdot\\ \cdot&1
\end{array}
\\\hline\\[-3mm]
\left[\begin{array}{cc}
b_{54}&b_{56}\\
b_{74}&b_{76}
\end{array}\right]
- 
\left[\begin{array}{c}
b_{51}\\ b_{71}
\end{array}\right]
b_{31}^{-1}
\left[\begin{array}{cc}
b_{34}& b_{36}
\end{array}\right]
& 
\begin{array}{cc}
\cdot & \cdot\\\cdot&\cdot
\end{array}
\end{array}\right]\\
&= 
\left[\begin{array}{cc}
0 & I_2\\
\sch(\Phi_1,\Phi_2) & 0 
\end{array}\right]
\enspace.
\end{align*}
This concludes the proof of Corollary~\ref{cor-l}.%\qed
\subsection{Singular examples}\label{sec-singular}
We now show how Theorem~\ref{th-1} allows
to solve singular cases in Lidski\u\i's theorem (Corollary~\ref{cor-l}),
and we also illustrate the limitations of Theorem~\ref{th-1}.
\begin{example}
Consider the following classical degenerate example,
taken from~\cite[Section~2.22]{wilkilson} and~\cite[Eqn~1.1]{mbo97}:
\def\ddiam{\mbox{\Large$\diamond$}\!\!\!\cdot}
\begin{align}
\sA_\eps=\sA_0+\eps b,
\;\mrm{where}\;
\sA_0= \Nil(1,1:3,2)
= \left[\begin{array}{ccc|cc} 
\cdot & 1& \cdot & \cdot&\cdot\\
\cdot & \cdot & 1& \cdot&\cdot\\
\odot & \cdot & \cdot& \dbox&\cdot\\\hline
\cdot & \cdot & \cdot& \cdot&1\\
\dbox & \cdot & \cdot& \dbox&\cdot
\end{array}
\right],\;\mrm{and}\; b\in \C^{n\times n} 
\enspace .\label{e-ex-2}
\end{align}
Recall that all the dots (whether they are
surrounded by boxes or circles, or not)
represent $0$. If the entry $b_{31}$
corresponding to the circled position in~\eqref{e-ex-2},
is zero, $\Phi_1$ is singular,
and we cannot apply Lidski\u\i's theorem (Corollary~\ref{cor-l}).
However, Theorem~\ref{th-1} can be applied. We can write 
$(\sA_\eps)_{ij} \simeq a_{ij} \eps^{A_{ij}}$,
with
\[
a= \left[\begin{array}{ccc|cc} 
b_{11} & 1& b_{13} & b_{14}&b_{15}\\
b_{21} & b_{22} & 1& b_{24}&b_{25}\\
0 & b_{32} & b_{33}& b_{34}&b_{35}\\\hline
b_{41} & b_{42} & b_{43}& b_{44}&1\\
b_{51} & b_{52} & b_{53}& b_{54}&b_{55}
\end{array}\right] \enspace,\;
\mrm{and}\;
A= 
\left[\begin{array}{ccc|cc} 
1 & 0& 1 & 1&1\\
1 & 1 & 0& 1&1\\
\infty & 1 & 1& 1&1\\\hline
1 & 1 & 1& 1&0\\
1& 1 & 1& 1&1
\end{array}\right] \enspace.
\]
We have $\alpha_1=\rhomin(A)=2/5$, $C_1=\{1,2,3,4,5\}$,
and since the critical graph of $A$, which is composed
only of the circuit $(1\to 2 \to 3 \to 4 \to 5 \to 1)$
covers all the nodes, we have $\GC(\ainf)=\GC(A)$. Thus,
for $G=\GC(\ainf)$,
\[
a^{G}= \left[\begin{array}{ccc|cc} 
\cdot & 1& \cdot & \cdot&\cdot\\
\cdot & \cdot & 1& \cdot & \cdot\\
\cdot & \cdot & \cdot& b_{45}& \cdot\\\hline
\cdot & \cdot & \cdot& \cdot&1\\
b_{51} & \cdot &  \cdot &  \cdot &  \cdot 
\end{array}\right] \enspace.
\]
Theorem~\ref{th-1} shows that,
if $b_{45}b_{51}\neq 0$, $\sA_\eps$
has five roots with asymptotics 
\[
\sL_\eps \sim \xi \eps^{2/5} \enspace, \; \mrm{where}\;
\xi^5= b_{45}b_{51} \enspace .
\]
The asymptotics of the eigenvectors can also
be obtained from Theorem~\ref{th-2}
(the computations are similar
to the case of Example~\ref{ex-1-eigenvect}).
\end{example}
\begin{example}
Let us discuss the following singular version of 
the illustrating example of~\cite{mbo97}.
Let $\sA_\eps = \sA_0 + \eps b$, where
$\sA_0 = \Nil(2,1,1:3,2,1)$, so that,
setting $G=G^c(\hat{A})$,
\[
a^{G}=
\def\mydot{\;\cdot\;}
\left[\begin{array}{ccc|ccc|cc|c}
\mydot & 1    & \mydot  & \mydot & \mydot & \mydot & \mydot&\mydot& \mydot\\
\mydot & \mydot& 1   & \mydot & \mydot & \mydot & \mydot&\mydot& \mydot\\
b_{31} & \mydot& \mydot   & b_{34} & \mydot & \mydot & b_{37}
&\mydot& b_{39}
\\\hline
\mydot & \mydot& \mydot   & \mydot & 1 & \mydot & \mydot&\mydot& \mydot\\
\mydot & \mydot& \mydot   & \mydot & \mydot & 1 & \mydot&\mydot& \mydot\\
b_{61} & \mydot& \mydot   & b_{64} & \mydot & \mydot & b_{67}&\mydot& b_{69}
\\\hline
\mydot & \mydot& \mydot   & \mydot & \mydot & \mydot & \mydot&1& \mydot\\
b_{81} & \mydot& \mydot   & b_{84} & \mydot & \mydot & b_{87}&\mydot& b_{89}
\\\hline
b_{91} & \mydot& \mydot   & b_{94} & \mydot & \mydot & b_{97}&\mydot& b_{99}
\end{array}
\right] \enspace .
\]
Consider the singular case where $b_{61}=b_{64}=0$.
We may keep $A$ as in Section~\ref{sec-cor},
but this gives little information
since $t^1$ is not invertible. 
However, $(\sA_{\eps})_{ij}\simeq a_{ij}\eps^{A_{ij}}$
still holds if we change the following values of $A$:
$A_{61}=A_{64}=\infty$. 
Then, we find 
$\alpha_1=1/3$, $C_1=\{1,2,3\}$, $\alpha_2=2/5$,
$C_2=\{4,5,6,7,8\}$, $\alpha_3=4/5$, $C_3=\{9\}$,
and the critical graphs $G^c_\ell (A)$, $\ell=1,2,3$ are represented
as follows:
\begin{center}
\begin{tabular}[c]{c}
\input{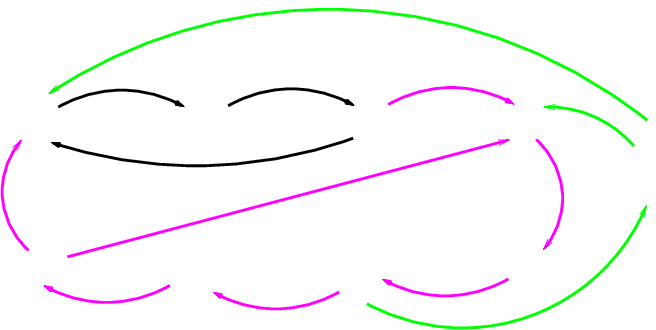}
\end{tabular}
\end{center}
with the same colouring convention as in Example~\ref{ex-critic}.

The matrix $t^1$ is invertible if, and only if, $b_{31}\neq 0$.
In this case, $\sA_{\eps}$ has three
eigenvalues with asymptotics
\(
\sL_{\eps}\sim\lambda \eps^{1/3}\),
corresponding to the different cubic roots $\lambda$ 
of $b_{31}$. We have 
\[
s^2=
\def\mydot{\cdot}
\left[\begin{array}{ccc|cc|c}
 \mydot & 1 & \mydot & \mydot&\mydot& \mydot\\
 \mydot & \mydot & 1 & \mydot&\mydot& \mydot\\
\mydot & \mydot & \mydot & b_{67}&\mydot& \mydot
\\\hline
 \mydot & \mydot & \mydot & \mydot&1& \mydot\\
 b'_{84} & \mydot & \mydot & \mydot&\mydot& b'_{89}
\\\hline
b'_{94} & \mydot & \mydot & \mydot&\mydot& b'_{99}
\end{array}
\right]
\]
where for instance
$b'_{84}=b_{84}-b_{81}b_{31}^{-1}b_{34}$.
Thus,  $t^2$ is invertible, if, and only if, $
b'_{84}b_{67}\neq 0$. When this is the case,
$\sA_{\eps}$
has five eigenvalues with asymptotics
\(
\sL_{\eps}\sim\lambda \eps^{2/5} 
\),
corresponding to the different quintic roots
$\lambda$ of $b'_{84}b_{67}$. 

The last critical graph, $G^c_3(A)$, that we just represented above,
does not have a disjoint circuit cover.
To see this, observe that there is no arc
from the set $\{3,8,9\}$ to the set $\{2,3,5,6,7,8,9\}$,
remark that the sum of the numbers of elements of these two sets,
which is $3+7=10$, exceeds the dimension of the matrix, which is $9$, and
apply the Frobenius-K\"onig theorem (see for instance~\cite[Th.~2.14]{br}).
Then, we know from Theorems~\ref{th-gam-bet} and~\ref{termrk}
that the greatest root, $\gamma_9$, of the min-plus characteristic
polynomial of $A$, $P_A$, is strictly greater than the greatest critical
value, $\beta_9=\alpha_3=4/5$, and by Theorem~\ref{th-gen=}, 
the exponent $\Lambda_9$ of the remaining eigenvalue of $\sA_\epsilon$
must be strictly greater than $4/5$. Thus, in this
case, Theorem~\ref{th-1} does not predict
the exponent of the eigenvalue of $\sA_\epsilon$ of minimal
modulus. However, this exponent can be obtained
as follows.  We already know that $\alpha_1=1/3$ and 
$\alpha_2=2/5$
are roots of respective multiplicity $3$ and $5$
of $P_A$, so the associated characteristic polynomial function
is of the form
$\widehat P_A(y)= (y\oplus \alpha_1)^3(y\oplus\alpha_2)^5(y\oplus \gamma_9)$.
One can check that $\perm A=4$, and since
$\hat P_A(\zero)=\perm A$, 
we deduce that $\alpha_1^3\alpha_2^5\gamma_9=3\gamma_9=4$,
therefore, $\gamma_9=1$. Then, one can derive from
Theorem~\ref{th-gen=} that $\Lambda_9=1$, for generic values
of $b$. The problem of finding the
leading coefficient of the corresponding
eigenvalue of $\sA_\epsilon$ is solved by
the result of~\cite{abg04}.
\end{example}
\begin{example}
Corollary~3.3 of~\cite{mbo97} identifies a special situation
where the leading exponent of a group of eigenvalues
can be found although the corresponding
matrix $\Phi_{\ell-1}$ appearing in Lidski\u\i's theorem
(see Corollary~\ref{cor-l} above) is not invertible.
We next give an example which cannot be solved
using the method of~\cite{mbo97} but which is
solved by Theorem~\ref{th-1}.
Let
\begin{align*}
\sA_\epsilon=
 \left[\begin{array}{cc|c|cc|c|c} 
\cdot & 1& \cdot & \cdot&\cdot&\cdot&\cdot\\
\cdot & \cdot & \epsilon b_{23}& \cdot&\cdot&\cdot&\cdot\\\hline
\epsilon b_{31} & \cdot & \cdot& \epsilon b_{34}&\cdot&\cdot&\cdot\\\hline
\cdot & \cdot & \cdot& \cdot&1&\cdot&\cdot\\
\cdot & \cdot & \cdot & \cdot & \cdot & \epsilon b_{56}&\cdot\\\hline
\cdot & \cdot & \cdot & \cdot & \cdot & \cdot&\epsilon b_{67}\\\hline
\epsilon b_{71} & \cdot & \cdot & \epsilon b_{74} & \cdot & \cdot&\cdot
\end{array}
\right] ,
\end{align*}
where the $b_{ij}$ are complex numbers, and, as above, the 
dots represent zero entries. The matrix $\sA_\epsilon$
is of the form $\sA_0+\epsilon b$, where the matrix $\sA_0$
is a nilpotent matrix conjugate to $\Nil(2,3:2,1)$.
The corresponding 
matrix $A$ is 
\[
A= \left[\begin{array}{cc|c|cc|c|c} 
\infty & 0& \infty & \infty & \infty & \infty & \infty  \\
\infty & \infty & 1& \infty& \infty& \infty& \infty\\\hline
1 & \infty & \infty & 1& \infty&\infty &\infty\\\hline
\infty & \infty & \infty& \infty&0&\infty&\infty\\
\infty & \infty & \infty & \infty & \infty & 1 & \infty\\\hline
\infty & \infty & \infty & \infty & \infty & \infty & 1\\\hline
1 & \infty & \infty & 1 & \infty & \infty & \infty
\end{array}\right] \enspace.
\]
We have $\alpha_1=2/3$, $C_1=\{1,2,3\}$, $\alpha_2=3/4$, $C_2=\{4,5,6,7\}$.
The critical graphs $G^c_1(A)$ and $G^c_2(A)$ are the following:
\begin{center}
\begin{tabular}[c]{c}
\input{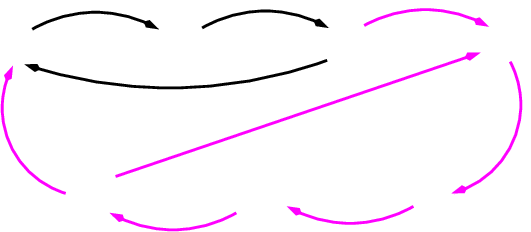}
\end{tabular}
\end{center}
with the same colouring convention as in the previous examples.
Theorem~\ref{th-1} shows that, for generic values of the coefficients $b_{ij}$,
the matrix $\sA_\eps$ has three eigenvalues
$\sL_\epsilon\sim \lambda\epsilon^{2/3}$, 
where $\lambda$ is a cubic root of $b_{23}b_{31}$,
and four eigenvalues
of the form $\lambda \epsilon^{3/4}$, where $\lambda$
is a quartic root of
\begin{align*}
b_{56}b_{67}(b_{74}-b_{71}b_{31}^{-1}b_{34}) \enspace .
\end{align*}
The following picture represents
the actual Newton polygon of the characteristic
polynomial $\det(\iY I-\sA_\epsilon)$,
for generic values of the $b_{ij}$ (this is
exactly the graph of $\divex{P}$, where $P=\perm(\iY I\oplus A)$.)
A monomial $\iY^i\epsilon^j$ is represented by the point
of coordinates $(i,j)$. 
Integer points are represented by small crosses. 
The actual Newton polygon
(black broken line) consists of two segments of respective slopes $-3/4$ and $-2/3$,
joining the three circles.
The approximation of the Newton polygon provided by Lidski\u\i's theorem
is given by the dashed broken line.
\begin{center}
\begin{tabular}[c]{c}
%auto-ignore
\begin{picture}(0,0)%
\includegraphics{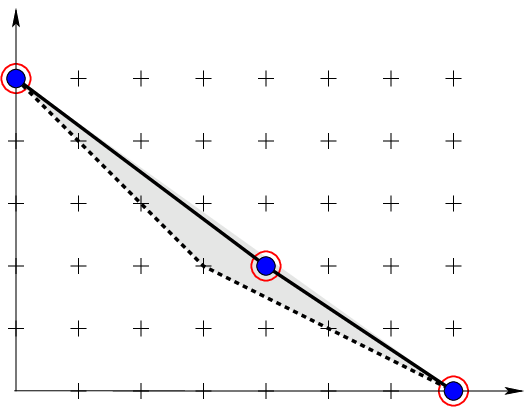}%
\end{picture}%
\setlength{\unitlength}{1973sp}%
\begingroup\makeatletter\ifx\SetFigFont\undefined%
\gdef\SetFigFont#1#2#3#4#5{%
  \reset@font\fontsize{#1}{#2pt}%
  \fontfamily{#3}\fontseries{#4}\fontshape{#5}%
  \selectfont}%
\fi\endgroup%
\begin{picture}(5042,3841)(7046,-5315)
\end{picture}%

\end{tabular}
\end{center}
The method of~\cite{mbo97} relies on the observation
that Lidski\u\i's theorem provides an approximation
of the Newton polygon, 
which is exact when the matrices $\Phi_\ell$ are invertible.
Corollary~3.3 of~\cite{mbo97} requires the absence
of integers points strictly between Lidski\u\i's approximation
and its chord, i.e.\ in the present case, in the interior
of the gray region.
Since this interior contains
the integer point $(4,2)$, the leading exponents $2/3$ and $3/4$
cannot be obtained from~\cite{mbo97}.
\end{example}

\def\cprime{$'$}

\end{document}